\documentclass[11pt]{article}
 \usepackage[margin=0.7in]{geometry} 
\usepackage{cite}
\usepackage{amsmath,amssymb,amsfonts}
\usepackage{algorithmic}
\usepackage{graphicx}
\usepackage{algorithm,algorithmic}
\usepackage{hyperref}
\hypersetup{hidelinks=true}
\usepackage{textcomp}
\usepackage{xz_style}

\renewcommand\footnotemark{}

\begin{document}
\title{Global Convergence of Receding-Horizon Policy Search \\
 in Learning Estimator Designs}
\author{Xiangyuan Zhang \quad Saviz Mowlavi \quad Mouhacine Benosman \quad Tamer Ba\c{s}ar 
\thanks{X. Zhang and T. Ba\c{s}ar are with the Department of ECE and CSL, University of Illinois Urbana–Champaign, Urbana, IL 61801 USA (e-mails: \{xz7, basar1\}@illinois.edu). S. Mowlavi and M. Benosman are with Mitsubishi Electric Research Laboratories, Cambridge, MA 02139 USA (e-mails: \{mowlavi, benosman\}@merl.com).}
\thanks{A preliminary version of this manuscript \cite{zhang2023learning} appeared in the proceedings of the 2023 American Control Conference and was presented on June 2, 2023, in San Diego, CA. }
}

\maketitle

\begin{abstract}
We introduce the receding-horizon policy gradient (RHPG) algorithm, the first PG algorithm with provable global convergence in learning the optimal linear estimator designs, i.e., the Kalman filter (KF). Notably, the RHPG algorithm does not require any prior knowledge of the system for initialization and does not require the target system to be open-loop stable. The key of RHPG is that we integrate vanilla PG (or any other policy search directions) into a dynamic programming outer loop, which iteratively decomposes the infinite-horizon KF problem that is constrained and non-convex in the policy parameter into a sequence of static estimation problems that are unconstrained and strongly-convex, thus enabling global convergence. We further provide fine-grained analyses of the optimization landscape under RHPG and detail the convergence and sample complexity guarantees of the algorithm. This work serves as an initial attempt to develop reinforcement learning algorithms specifically for control applications with performance guarantees by utilizing classic control theory in both algorithmic design and theoretical analyses. Lastly, we validate our theories by deploying the RHPG algorithm to learn the Kalman filter design of a large-scale convection-diffusion model. We open-source the code repository at \url{https://github.com/xiangyuan-zhang/LearningKF}.
\end{abstract}

\section{Introduction}
In recent years, policy-based reinforcement learning (RL) methods \cite{sutton2000policy, kakade2002natural,schulman2015trust,schulman2017proximal} have gained increasing attention in continuous control applications \cite{schulman2015high,lillicrap2015continuous,recht2019tour}. While traditional model-based techniques synthesize controller designs in a case-by-case manner \cite{anderson1979optimal, anderson1990optimal}, model-free policy gradient (PG) methods promise a universal framework that learns controller designs in an end-to-end fashion. The universality of model-free PG methods makes them desired candidates in complex control applications that involve nonlinear system dynamics and imperfect state measurements. Despite countless empirical successes, the theoretical properties of model-free PG methods still need to be thoroughly investigated in continuous control. Initiated by \cite{fazel2018global}, a recent line of research has well-analyzed the sample complexity of zeroth-order PG methods in several linear state-feedback control benchmarks, including linear-quadratic regulator (LQR) \cite{fazel2018global, mohammadi2019convergence, malik2020derivative, hambly2020policy, perdomo2021stabilizing, ju2022model, zhang2023revisiting}, distributed/decentralized LQR \cite{li2019distributed, furieri2019learning}, and linear robust control \cite{gravell2019learning, zhang2019policymixed, zhang2021derivative}. However, the theoretical properties of PG methods remain elusive in the output-feedback control settings, where the state measurement process could be corrupted by statistical noises and/or other (possibly adversarial) disturbances. 

In this work, we study the convergence and sample complexity of PG methods in the discrete-time infinite-horizon Kalman filtering (KF) problem \cite{kalman1960new, anderson1979optimal}. Recognized as one of the cornerstones of modern control theory \cite{basar2001control}, the KF problem aims to generate optimal estimates of the unknown system states over time by utilizing a sequence of observed measurements corrupted by statistical noises. Furthermore, in the linear-quadratic Gaussian (LQG) problem, the separation principle \cite{astrom1971introduction} states that the optimal control law combines KF and LQR. Thus, KF is a fundamental benchmark for studying the sample complexity of model-free PG methods beyond state-feedback settings. 

Despite being the dual problem to noise-less LQR \cite{kalman1960general, astrom1971introduction}, the KF problem possesses a substantially more complicated optimization landscape from the model-free PG perspective since the KF itself is a dynamical system rather than a static matrix. Specifically, the optimization problem over dynamic filters might admit multiple suboptimal stationary points, and the optimal KF possesses a set of equivalent realizations up to similarity transformations \cite{zheng2021analysis, umenberger2022globally}. None of the above challenges appear when using model-free PG to learn a static LQR policy \cite{fazel2018global, mohammadi2019convergence, malik2020derivative, hambly2020policy, perdomo2021stabilizing, ju2022model, zhang2023revisiting}. As a result of the challenging landscape the filtering problem presents, only a few papers have focused on dynamic output-feedback settings. In particular, \cite{zheng2021analysis} has analyzed the optimization landscape of LQG, and \cite{umenberger2022globally} has shown that an informativity-regularized PG method provably converges to an optimal filter in the continuous-time KF problem, assuming that the model is known. However, \cite{umenberger2022globally} has assumed that the target system is open-loop stable and assumed that the control engineer has prior knowledge of a filter that satisfies an informativity condition. It is also unclear if the techniques in \cite{umenberger2022globally} can be directly applied to the model-free setting and result in any sample complexity guarantees. Thus, obtaining sample complexity of model-free PG methods in the KF problem has remained a significant challenge.

This work addresses these challenges by introducing a receding-horizon PG  (RHPG) algorithm and establishing its global convergence and sample complexity. In contrast to direct policy search, the RHPG algorithm integrates vanilla PG (or any other policy search directions) into a dynamic programming (DP) outer loop, which iteratively decomposes the infinite-horizon KF problem that is constrained and non-convex in the policy parameter into a sequence of static estimation problems that are unconstrained and strongly-convex. Then, we show that solving the sequence of static estimation problems results in the global convergence of RHPG toward the KF, which is the optimal linear filter. We further establish the total sample complexity of the RHPG to be $\tilde{\mathcal{O}}(\epsilon^{-2})$ for the learned filter to be $\epsilon$-close in policy distance to KF, which is the first sample complexity result of PG methods in the output-feedback control settings. Notably, the RHPG algorithm does not require any prior knowledge of the system to generate a valid initialization and does not require the target system to be open-loop stable. This removes two restrictive assumptions in the previous work \cite{umenberger2022globally}. We validate our theories by learning the KF design of a large-scale convection-diffusion model.

Compared to the preliminary results included in \cite{zhang2023learning}, this work presents a comprehensive study of RHPG as a model-free RL approach in learning estimator designs. In particular, our contributions, in addition to those listed in \cite{zhang2023learning}, are three-folded. First, we analyze the optimization landscape in Theorem \ref{theorem:quadratic}, which clarifies the properties of quadratic programs in RHPG and provides a theoretical foundation for selecting the algorithmic parameters of PG methods. Second, we discuss the insights of the RHPG design in Sec. \ref{sec:alg_design} and compare RHPG with standard PG methods regarding parametrization and landscape, computational efficiencies, and requirements on the simulation oracles. These discussions provide the intuitions behind the mathematical developments of RHPG and should benefit future research on improving the algorithm design and the theoretical analysis. Lastly, we open-source numerical experiments on the data-driven estimation of a large-scale dynamical system, which corroborates the theories and demonstrates both the effectiveness and scalability of RHPG.

This work attempts to develop RL algorithms specifically for control and estimation tasks with performance guarantees, by utilizing classic control theory in both algorithmic design and theoretical analyses. The dual theories and implementations of RHPG to the LQR problem have been presented in \cite{zhang2023revisiting}. Through this line of work, we demonstrate the significant utilization of DP in overcoming the challenging optimization landscape and streamlining the analyses when deploying model-free PG methods to linear control and estimation tasks. Due to the separation principle \cite{astrom1971introduction}, our results shed light on applying model-free PG methods in solving the LQG problem through a sequential design of controller and estimator.

The structure of the paper is as follows. In Sec. \ref{sec:formulation}, we define the infinite- and finite-horizon settings of the KF problem and formulate them as policy optimization problems. In Sec. \ref{sec:alg}, we introduce the RHPG algorithm and provide the general theory and intuition that backs its design. In Sec. \ref{sec:landscape}, we analyze the optimization landscape of solving the KF problem using RHPG and establish the global convergence and sample complexity guarantees of the algorithm. Lastly, we present the numerical studies on a convection-diffusion model in Sec. \ref{sec:sim}. The paper ends with the concluding remarks of Sec. \ref{sec:conclude}, and an appendix that contains formal proofs of the main results.

\subsection{Notations}\label{sec:notations}
For a square matrix $X$, we denote its trace, spectral norm, condition number, and spectral radius by $\Tr(X)$, $\|X\|$, $\kappa_X$, and $\rho(X)$ resp. We define the $W$-induced norm of $X$ as $\|X\|^2_W := \max_{z \neq 0} \frac{z^{\top}X^{\top}WXz}{z^{\top}Wz}$. If $X$ is further symmetric,  we use $X > 0$, $X\geq 0$, $X\leq 0$, and $X <0$ to denote that $X$ is positive definite (pd), positive semi-definite (psd), negative semi-definite (nsd), and negative definite (nd), resp. We use $x\sim \cN(\mu,\Sigma)$ to denote a Gaussian random vector with mean $\mu$ and covariance $\Sigma$. Lastly, we use $\bI$ and $\bm{0}$ to denote the identity and zero matrices, resp., with appropriate dimensions.

\section{Preliminaries}\label{sec:formulation}
\subsection{Infinite-Horizon Kalman Filtering}
Consider the discrete-time linear time-invariant system
\begin{align}\label{eqn:inf_LQE_dynamics}
	x_{t+1} = Ax_t+ w_t, \quad y_t = Cx_t + v_t,
\end{align}
where $x_t \in \RR^n$ is the state, $y_t \in \RR^m$ is the output measurement, and $w_t \sim \cN(\bm{0}, W)$, $v_t\sim \cN(\bm{0}, V)$ are sequences of i.i.d. zero-mean Gaussian noises for some $W, V > 0$, also independent of each other. The initial state is also assumed to be a Gaussian random vector such that $x_0 \sim \cN(\bar{x}_0, X_0)$, independent of $\{w_t, v_t\}$, with $\bar{x}_0 \neq \bm{0}$ and $X_0 > 0$. Additionally, we assume that $(C, A)$ is observable and note that the condition $W>0$ readily leads to controllability of $(A, W^{1/2})$, which is a standard condition in KF. 

The KF problem aims to generate a sequence of estimated states, denoted by $\hat{x}_t$ for each $t$, that minimizes the infinite-horizon mean-square error (MSE):
\begin{align}\label{eqn:inf_mse}
	\cJ_{\infty} := \lim_{N\to \infty}\frac{1}{N}\EE\bigg\{\sum_{t=0}^{N}(x_t-\hat{x}_t)^{\top}(x_t-\hat{x}_t)\bigg\}.
\end{align}
Moreover, each $\hat{x}_t$ can only depend on the history and output measurements up to but not including $t$, i.e., $\{y_0, \cdots, y_{t-1}\}$. The celebrated result of Kalman \cite{kalman1960new} showed that the $\cJ_{\infty}$-minimizing filter (could also be called $1$-step predictor), which exists under the controllability and the observability conditions, has the form of 
\begin{align}\label{eqn:dynamic_filter}
	\hat{x}_{t+1}^* &= (A-L^*C)\hat{x}_t^* + L^*y_t, \quad  \hat{x}_{0}^* = \bar{x}_0, \\ 
	L^*&=A\Sigma^* C^{\top}(V + C\Sigma^* C^{\top})^{-1}, \label{eqn:kalman_gain}
\end{align} 
where $L^*$ is the Kalman gain and $\Sigma^*$ represents the unique pd solution to the filter algebraic Riccati equation (FARE):
\begin{align}
	\Sigma &= A\Sigma A^{\top} - A\Sigma C^{\top}(V + C\Sigma C^{\top})^{-1}C\Sigma A^{\top} + W. \label{eqn:filter_riccati}
\end{align}
Hence, without any loss of optimality, we can restrict the search to the class of filters of the form $\hat{x}_{t+1} = A_L\hat{x}_t + B_Ly_t$ and then parametrize the KF problem as a minimization problem over $A_L$ and $B_L$ subject to a stability constraint\footnote{Extending the results in this work to the setting with instantaneous feedback measurement (i.e., allowing $\hat{x}_t$ to depend also on $y_t$, and hence replacing $y_t$ in (\ref{eqn:filter_form}) with $y_{t+1}$) would be straightforward.}
\begin{align}\label{eqn:filter_form}
	&\min_{A_L, B_L} \quad \cJ_{\infty}(A_L, B_L) \quad \text{s.t.} \quad \hat{x}_{t+1} = A_L\hat{x}_t + B_Ly_t ~\ \text{and} ~\ \rho(A_L) < 1.
\end{align}
Note that by \eqref{eqn:dynamic_filter}, there indeed exists a solution to \eqref{eqn:filter_form} where $(A_L^*, B_L^*)=(A-L^*C, L^*)$. Note also that when the pair $(A, C)$ is known, \eqref{eqn:filter_form} involves an over-parametrization since solving \eqref{eqn:filter_form} is equivalent to optimizing a single variable $B_L$. However, in the model-free setting where $(A, C)$ is unknown, which is the target setting of our paper, it is reasonable to parametrize the KF problem as in \eqref{eqn:filter_form}. Until now, obtaining sample complexity of model-free PG methods in solving the KF problem \eqref{eqn:filter_form} has remained a major challenge.

\begin{figure*}[t]
\centering 
\includegraphics[width = 0.85\textwidth]{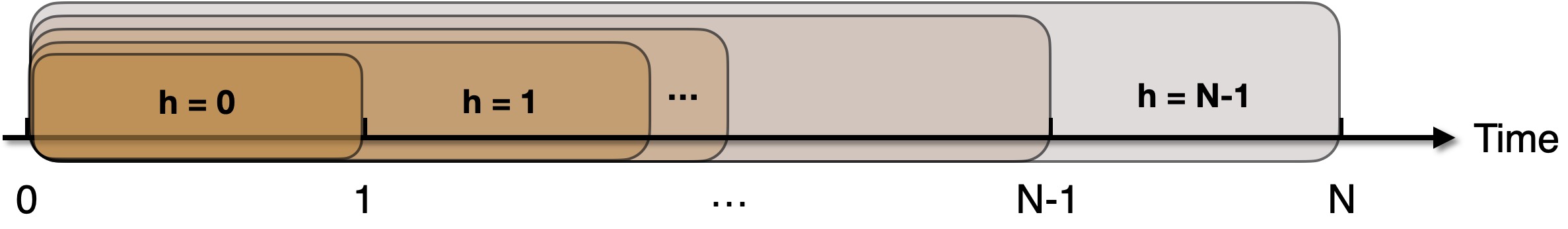}	
\caption{Algorithm \ref{alg:DPfilter}, executed forward in time, constructs an $(h+1)$-horizon KF problem from $t=0$ to $t=h+1$ at each iteration indexed by $h$.}\label{fig:DP_illustration}
\end{figure*}

\subsection{Finite-Horizon Kalman Filtering}
We now discuss the finite-$N$-horizon KF problem, also described by the system dynamics \eqref{eqn:inf_LQE_dynamics}. Adopting the same parametrization as in \eqref{eqn:filter_form}, but this time allowing time-dependence, and again without any loss of optimality, we represent the finite-horizon KF problem as a minimization problem over a sequence of time-varying filter parameters $\{A_{L_t}, B_{L_t}\}$, for all $t\in\{0, \cdots, N-1\}$, \par 
\small \vspace{-1em}
\begin{align}\label{eqn:finite-mse}
	\min_{\{A_{L_t}, B_{L_t}\}}  & \cJ\big(\hspace{-0.1em}\{A_{L_t}, B_{L_t}\}\hspace{-0.1em}\big) \hspace{-0.2em}:=\hspace{-0.2em} \EE\bigg\{\hspace{-0.15em}\sum_{t=0}^{N}(x_t-\hat{x}_t)^{\hspace{-0.1em}\top}\hspace{-0.1em}(x_t-\hat{x}_t)\hspace{-0.15em}\bigg\} \\
	&\text{s.t.} \quad \hat{x}_{t+1} = A_{L_t}\hat{x}_t + B_{L_t}y_t, \quad \hat{x}_0 = \bar{x}_0. \nonumber
\end{align}
\normalsize
The minimum in \eqref{eqn:finite-mse} can be achieved by $(A^*_{L_t}, B^*_{L_t}) = (A-L^*_tC, L^*_t)$, where $L^*_t$ is the time-varying Kalman gain \par 
\small \vspace{-1em}
\begin{align}
	&L_t^*=A\Sigma_t^* C^{\top}(V + C\Sigma_t^* C^{\top})^{-1}, \quad  \Sigma_0^* = X_0.\label{eqn:kalman_gain_finite}\\
	\hspace{-0.2em} \Sigma^*_{t+1} &= A\Sigma^*_{t}A^{\top} \hspace{-0.1em}-\hspace{-0.1em} A\Sigma^*_t C^{\top}(V \hspace{-0.1em}+\hspace{-0.1em} C\Sigma^*_t C^{\top})^{-1}C\Sigma^*_t A^{\top} \hspace{-0.4em}+\hspace{-0.1em} W. \label{eqn:filter_riccati_finite}
\end{align}
\normalsize
The solutions $\Sigma^*_t$, for all $t \in \{0, \cdots, N-1\}$, generated by the filter Riccati difference equation (FRDE) \eqref{eqn:filter_riccati_finite} always exist and are unique and pd, due to $V > 0$, $W > 0$, and the iteration starts with $\Sigma^*_0 = X_0 > 0$.

\section{Receding-Horizon Policy Gradient}\label{sec:alg}
\subsection{Kalman Filtering and Dynamic Programming}
It is well known that the solution of the FRDE \eqref{eqn:filter_riccati_finite} converges monotonically to the stabilizing solution of the FARE \eqref{eqn:filter_riccati} at an exponential rate \cite{chan1984convergence, hassibi1999indefinite}. Then, it readily follows that the optimal time-varying filter $(A_{L_t}^*, B_{L_t}^*)$ to the finite-horizon KF problem \eqref{eqn:finite-mse} also converge monotonically to the time-invariant $(A_{L}^*, B_{L}^*)$ as $N \to \infty$. We present this convergence result in the following theorem, which plays an important role in our algorithm design.

\begin{theorem}\label{lemma:finite_approximation}
The finite-horizon Kalman gain as in \eqref{eqn:kalman_gain_finite} converges  to the infinite-horizon Kalman gain defined in \eqref{eqn:kalman_gain} exponentially fast as $N\to \infty$. Specifically, using $\|\cdot\|_*$ to denote the $\Sigma^*$-induced norm and letting
	\begin{align}\label{eqn:N0}
		N_0 = \frac{1}{2}\cdot \frac{\log\big(\frac{\|X_0-\Sigma^*\|_*\cdot\kappa_{\Sigma^*}\cdot \|A_L^*\|\cdot\|C\|} {\epsilon\cdot\lambda_{\min}(V)}\big)}{\log\big(\frac{1}{\|A_L^*\|_*}\big)} + 1.
	\end{align}
	where $\|A_L^*\|_* <1$, we have that, for all $N\geq N_0$, it holds that $\|L^*_{N-1} - L^*\| \leq \epsilon$ for any $\epsilon > 0$. If additionally $X_0 > \Sigma$ holds, then $L^*_t$ is stabilizing for all $t \geq 0$ in the sense that $\rho(A_{L^*_t}) = \rho(A-L^*_tC) < 1$. 
\end{theorem}\par  
The proof of Theorem \ref{lemma:finite_approximation} is provided in \S\ref{proof:finite}. Theorem \ref{lemma:finite_approximation} quantifies how different system parameters affect the non-asymptotic convergence rate of the time-varying filters to the time-invariant KF. It further demonstrates that if $N\sim \cO(\log(\epsilon^{-1}))$, then the filter $(A_{L_{N-1}}^*, B_{L_{N-1}})$ will be $\epsilon$-close to the infinite-horizon KF ($A_L^*, B_L^*$). Furthermore, if $\epsilon$ is sufficiently small (i.e., smaller than the stability margin of $A_L^*$), then it holds that $\rho(A_{L_{N-1}}^*)<1$. When $X_0 > \Sigma$, the system is sufficiently excited, and as a result, the frozen filter at any $t \geq 0$ is stable in the sense that $\rho(A_{L^*_t}) < 1$.

\subsection{Algorithm Design}\label{sec:alg_design}
Instead of solving the infinite-horizon KF problem \eqref{eqn:finite-mse} directly, we introduce the RHPG algorithm, which first selects a sufficiently large problem horizon $N$ according to Theorem \ref{lemma:finite_approximation}, then constructs and solves $N$ static estimation problems sequentially (see Figure \ref{fig:DP_illustration} for an illustration) using PG methods. We describe the procedure of the RHPG algorithm below. 

\renewcommand\arraystretch{1.15}
\begin{algorithm}[H]
\caption{Receding-Horizon Policy Gradient (RHPG)}\label{alg:DPfilter}
\begin{algorithmic}[1]
\renewcommand{\algorithmicrequire}{\textbf{Input:}}
 \renewcommand{\algorithmicensure}{\textbf{Output:}}
 \REQUIRE Problem horizon $N$
 \STATE Initialize $A_{L_{0}},  B_{L_{0}} \leftarrow \bm{0}_{n\times n}, \bm{0}_{n\times m}$
\FOR{$h \in \{0, \cdots, N-1\}$}
        \STATE Solve \eqref{eqn:induction} using PG methods until convergence
        	\STATE Use the convergent filter $A_{L_h}, B_{L_h}$ to warm-start PG updates for the next iteration
        	\ENDFOR

\STATE \hspace{0.5cm}\textbf{return} $A_{L_{N-1}},  B_{L_{N-1}}$
\end{algorithmic}
\label{alg1}
\end{algorithm}

We provide detailed implementations of zeroth-order and first-order RHPG in Algorithms \ref{alg:zeroth_order_GD} and \ref{alg:analytic_gd}, respectively. 

The RHPG algorithm is executed forward in time, and in the first iteration, the algorithm learns the optimal filter for a one-step static estimation problem. Then, every subsequent RHPG iteration extends the problem horizon by adding one additional time step after the initial one. Formally, at each iteration indexed by $h$, the RHPG algorithm constructs an $(h+1)$-horizon KF problem from $t=0$ to $t=h+1$, but we fix the filter parameters for all $t\in\{0, \cdots, h-1\}$ as those generated from earlier iterations and only optimize for latest filter parameters $(A_{L_h}, B_{L_h})$. This renders each iteration of the RHPG algorithm into solving a static estimation problem that is quadratic in $(A_{L_h}, B_{L_h})$.

Mathematically, for every $h\in\{0, \cdots, N-1\}$, the RHPG algorithm solves the following minimization problems
\begin{align}
	&\hspace{-1.7em}\min_{A_{L_h}, B_{L_{h}}} \hspace{-0.1em}\cJ_h \hspace{-0.1em}:=\hspace{-0.1em} \EE_{x_0, w_t, v_t, \theta_0} \Big\{\hspace{-0.1em} \sum_{t=0}^{h+1}(x_t-\hat{x}_t)^{\hspace{-0.1em}\top}\hspace{-0.1em}(x_t-\hat{x}_t)\hspace{-0.1em}\Big\} \label{eqn:induction}\\
	\text{s.t.} ~ &\hat{x}_{t+1} = A_{L_t}^*\hat{x}_t + B_{L_t}^*y_t, \ \forall t \in \{0, \cdots, h-2\}, \ \hat{x}_0 = \bar{x}_0\nonumber\\
	&x_{t+1} = Ax_t + w_t, \ \forall t \in \{0, \cdots, h-2\}, \ x_0 \sim \cN(\bar{x}_0, X_0)\nonumber\\
	 &\hat{x}_h =  A_{L_{h-1}}^*\hat{x}_{h-1} + B_{L_{h-1}}^*y_{h-1} + \theta_0 \label{eqn:inject1}\\
	 &x_{h} = Ax_{h-1} + w_{h-1} + \theta_0 \label{eqn:inject2}.
\end{align}
where $\theta_0 \sim \cN(\bm{0}, \Theta) \in \RR^{n}$ is sampled independently to $x_0, w_t, v_t$ and satisfies $\Theta > 0$. The purpose of injecting an additional ``small'' noise $\theta_0$ in \eqref{eqn:inject1}-\eqref{eqn:inject2} is to ensure the strict convexity of the quadratic program \eqref{eqn:induction}; we will formally justify it in Sec. \ref{sec:landscape}. 

\begin{figure}
	\centering
\includegraphics[width=0.55\textwidth]{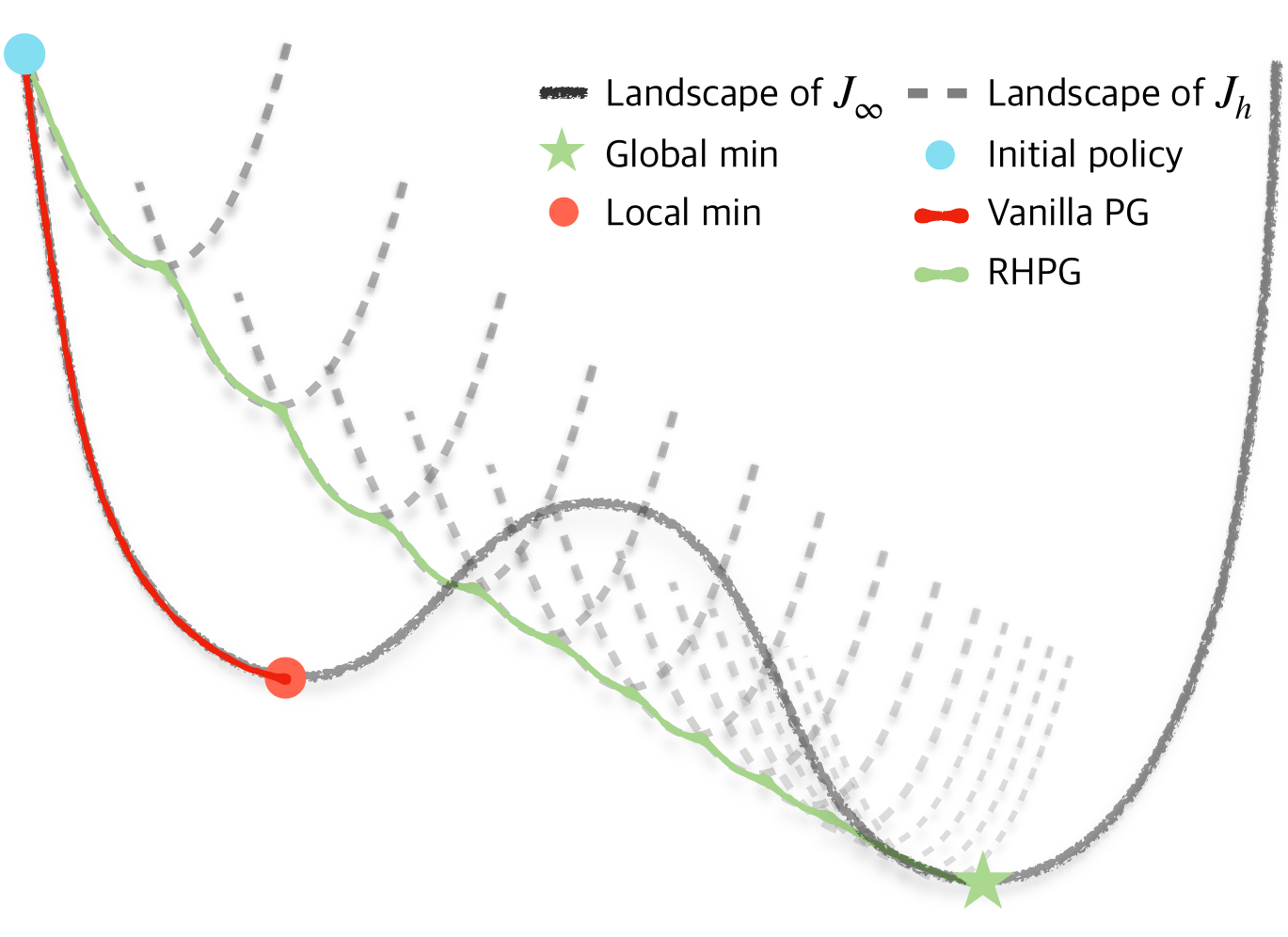}
	\caption{Illustration of the different optimization landscapes and search directions under vanilla PG and RHPG. *Note that this graph does not imply that the cost values of $\cJ_h$ are higher than that of $\cJ_{\infty}$. }\label{fig:landscape}
\end{figure}

Before presenting the theoretical analyses, we provide a few comments and discussions regarding the algorithmic designs of RHPG. From the optimization landscape perspective and compared to the vanilla PG update that is trapped by suboptimal first-order stationary points, RHPG provides a new policy search direction that points toward the global optimum by learning the KF design step-by-step, where each subproblem admits a unique global optimum due to the quadratic landscape. See Figure \ref{fig:landscape} for an illustration. Our algorithmic design shares a similar flavor as the curriculum learning literature \cite{bengio2009curriculum, wang2021survey, soviany2022curriculum}, where an agent (controller/estimator) evolves by mastering the simplest tasks first and then gradually conquering tasks that are more and more challenging. In our case, the RHPG algorithm first learns the filter capable of predicting only the immediate next state from scratch but then keeps adapting/evolving to handle new filtering tasks with longer and longer problem horizons (cf., the warm-start step in Line $4$ of Algorithm \ref{alg:DPfilter}). When the problem horizon becomes sufficiently large (as characterized by Theorem \ref{lemma:finite_approximation}), the filter converges globally, and the behaviors of the infinite-horizon KF, such as closed-loop stability, begin to emerge. Lastly, since RHPG starts by learning the simplest static estimation task and every subproblem is unconstrained, RHPG does not require any specific filter for initializations. In other words, it suffices to initialize arbitrarily for searching a static estimator. 

On the computational side, it may seem at first glance that RHPG is less efficient compared to vanilla PG since it solves $N$ optimization problems instead of $1$. This, however, turns out not to be true. When applying (sampled-based) vanilla PG to the infinite-horizon objective $\cJ_{\infty}$ directly, the rollout length is typically a very large finite number so that an accurate PG estimate can be obtained. The rollout length for RHPG is $1$ in its first iteration, which coincides with the iteration that needs the largest number of PG steps since it learns an optimal one-step static estimator \emph{from scratch}. The rollout length remains very small in the first few iterations, where the time-varying filters \eqref{eqn:kalman_gain_finite} are distinct across time. When the rollout length becomes moderate, only a few PG updates are needed to fine-tune the filter; see Figure \ref{fig:landscape}. Due to the much shorter rollout trajectories, the computational efficiency of RHPG is comparable, if not better, to vanilla PG. We provide the complexity analysis in Sec. \ref{sec:landscape}.

Lastly, we discuss the requirements for the simulation oracle. To sample the gradients, we require the standard assumption that the user has access to a simulator such that for any input filter $(A_{L_h}, B_{L_h})$, the simulator can return an empirical value of the objective function \eqref{eqn:induction}. This requires the simulator to generate exact state trajectories of the simulated model, but it only reveals a noisy scalar objective value to the learning algorithm. The requirement is reasonable for the offline learning setting since the algorithm does not use any system information directly. However, building a simulator naturally requires knowledge of the system model, which could be exact, approximate, or simplified. Transferring the simulated policies to a real system (a.k.a., Sim2Real) that might exhibit different dynamics requires a \emph{provable} robustness guarantee of the learned controller/estimator, which is beyond the scope of the present paper, but it is an important current and future research topic \cite{pinto2017robust, zhang2019policymixed, zhang2020rarl, zhang2021derivative, zhang2023ifac, cui2023reinforcement}. 

\subsection{Bias of Model-Free Receding-Horizon Filtering}\label{sec:bias}

\begin{figure}[t]
	\centering
	\includegraphics[width=0.65\textwidth]{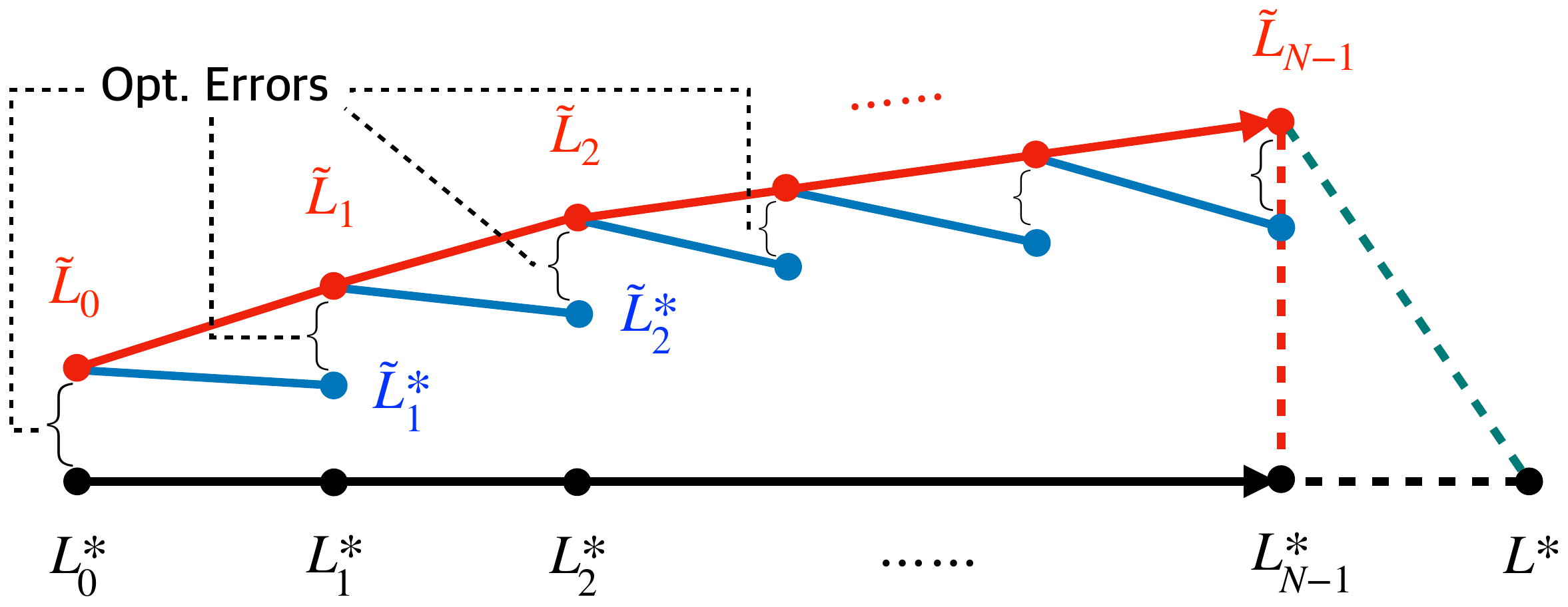}
	\caption{As illustrated, Theorem \ref{lemma:finite_approximation} bounds the distance between $L^*_{N-1}$ and $L^*$ and Theorem \ref{theorem:KF_DP} analyzes the forward propagation of the optimization errors from each iteration of RHPG. Combining two theorems upper-bounds the total policy gap to the infinite-horizon KF. }\label{fig:proof_sketch}
\end{figure}

The RHPG algorithm is backed by Bellman's principle of optimality, which requires solving each iteration exactly. However, iterative algorithms such as PG can only return an $\epsilon$-accurate solution in a finite time. We generalize the dynamic programming principle in the following theorem to analyze how computational errors accumulate in the forward DP process and provide an optimality guarantee for the filter that the RHPG algorithm returns. \par 

\begin{theorem}\label{theorem:KF_DP}
	Choose the problem horizon $N$ following Theorem \ref{lemma:finite_approximation} and assume that one can compute, for all $h\in\{0, \cdots, N-1\}$ and some $\epsilon > 0$, filter $(\tilde{A}_{L_h}, \tilde{B}_{L_h})$ that satisfies \par
	\small \vspace{-1em}
	\begin{align*}
		\big\|\tilde{A}_{L_h}  \hspace{-0.1em}-\hspace{-0.1em} \tilde{A}_{L_h}^*\big\|, \big\|\tilde{B}_{L_h} \hspace{-0.1em} -\hspace{-0.1em} \tilde{B}_{L_h}^*\big\| \hspace{-0.1em}\sim\hspace{-0.1em}\cO(\epsilon\cdot \texttt{poly}(\text{system parameters})),
	\end{align*}
	\normalsize
	where $(\tilde{A}_{L_h}^*, \tilde{B}_{L_h}^*)$ is the unique minimizing solution of $\cJ_h$ in \eqref{eqn:induction}, after setting the filters for all $t \in \{0, \cdots, h-2\}$ to be those computed in the previous iterations and are $\epsilon$-close to the minimum of $\cJ_0, \cdots, \cJ_{h-1}$, respectively. Then, the RHPG algorithm outputs $(\tilde{A}_{L_{N-1}}, \tilde{B}_{L_{N-1}})$ that satisfies $\big\|[\tilde{A}_{L_{N-1}} \ \tilde{B}_{L_{N-1}}] - [A^*_L \ B^*_L]\big\| \leq \epsilon$, where $(A^*_L, B^*_L)$ represents the infinite-horizon KF. If further $\epsilon$ is sufficiently small such that $\epsilon < 1-\|A_L^*\|_*$, then $\tilde{A}_{L_{N-1}}$ satisfies $\rho(\tilde{A}_{L_{N-1}}) < 1$.
\end{theorem}

We illustrate Theorem \ref{theorem:KF_DP} in Figure \ref{fig:proof_sketch} and defer its proof to \S\ref{proof:KF_DP}. Theorem \ref{theorem:KF_DP} guarantees that if every iteration of the DP is solved to an $\cO(\epsilon)$-accuracy, then the convergent filter after completing the $N$-step DP procedure is at most $\epsilon$-away from the exact infinite-horizon KF.  We note that the RHPG algorithm utilizes two layers of approximations. First, the solution to the infinite-horizon KF problem \eqref{eqn:filter_form} is approximated by the solution to a finite-horizon KF problem \eqref{eqn:finite-mse}, where we choose $N\sim \cO(\log(\epsilon^{-1}))$ due to the exponential attraction of the Riccati equation. Then, we solve the finite-horizon KF problem \eqref{eqn:finite-mse} by integrating forward DP with model-free policy search. Combining the two steps addresses the infinite-horizon KF task with a provable global convergence guarantee using only samples of system trajectories.

\section{Optimization Landscape, Convergence, and Sample Complexity}\label{sec:landscape}
We first present the optimization landscape of the static estimation problem \eqref{eqn:induction} in the following theorem. 

\begin{theorem}\label{theorem:quadratic}
	For every $h\in \{0, \cdots, N-1\}$, the quadratic objective $\cJ_h$ defined in \eqref{eqn:induction} is twice continuously differentiable, and its Hessian matrix can be represented as
	\begin{align*}
		H_{h} &= \EE_{x_h, \hat{x}_h, v_h}\begin{bmatrix}
		\hat{x}_{h}\hat{x}_{h}^{\top} & \hat{x}_{h}y_{h}^{\top} \\ y_{h}\hat{x}_{h}^{\top} & y_{h}y_{h}^{\top}
	\end{bmatrix} \\
	&= \begin{bmatrix}
		\mu_{\hat{x}_{h}}\mu_{\hat{x}_{h}}^{\top} +\Theta & (\mu_{\hat{x}_{h}}\mu_{x_{h}}^{\top} +\Theta)C^{\top} \vspace{0.2em}\\ C(\mu_{x_{h}}\mu_{\hat{x}_{h}}^{\top} +\Theta) & C(\mu_{x_{h}}\mu_{x_{h}}^{\top} +\Theta)C^{\top} + V
	\end{bmatrix} > 0,
	\end{align*}
	where we have used $\mu_{\hat{x}_{h}}$ and $\mu_{x_h}$ to denote $\EE[\hat{x}_h]$ and $\EE[x_h]$, respectively, and $\Theta > 0$ is the covariance matrix of the zero-mean Gaussian random vector $\theta_0$ in \eqref{eqn:inject1}-\eqref{eqn:inject2}. Moreover, the objective function $\cJ_h$ is strongly convex with constant $\lambda_{\min}(H_h)$ and smooth with constant $\lambda_{\max}(H_h)$. Lastly, introducing the additional random vector $\theta_0$ is without any loss of optimality in the sense that the time-varying KF characterized by \eqref{eqn:kalman_gain_finite}-\eqref{eqn:filter_riccati_finite} represents the unique minimum of \eqref{eqn:induction}. 
\end{theorem}

The proof of Theorem \ref{theorem:quadratic} is deferred to Sec. \ref{proof:quadratic}, where an extended discussion on the effect of $\theta_0$ is also provided. In short, introducing an additional  ``small''  $\theta_0$ in \eqref{eqn:inject1}-\eqref{eqn:inject2} ensures the strict convexity of the quadratic objective $\cJ_h$ with respect to $A_{L_h}$, while $\cJ_h$ is strictly convex in $B_{L_h}$ with or without $\theta_0$ due to the condition $V >0$. 

Denote $\pi_t = \big[A_{L_t} \mid B_{L_t}\big]$. We define the analytic vanilla PG of $\cJ_h$ for every $h \in \{0, \cdots, N-1\}$ to be\footnote{Note that one can use any other policy search directions such as natural PG \cite{kakade2002natural} or least-squares policy iteration \cite{lagoudakis2003least} to replace vanilla PG in RHPG.}
\begin{align}
	&\nabla_{\pi_h}\cJ_h(\pi_h) = 2 \Big[\pi_h(\Psi_h + \Delta) - (G_h + \Xi)\Big], \label{eqn:PG_RHPG_KF}
\end{align}
where
\begin{align}
&\Delta =
 \bigg[
\begin{array}{c|c}
  \Theta & \Theta C^{\top} \\
  \hline
  C\Theta & C\Theta C^{\top}
\end{array}
\bigg], \ \Xi = \big[\begin{array}{c|c}
  A\Theta & A\Theta C^{\top}
\end{array}\big] \label{eqn:KF_grad1}\\
& \Psi_t =  \bigg[
\begin{array}{c|c}
  \Var(\hat{x}_t) & \Cov(x_t, \hat{x}_t)^{\top} C^{\top} \\
  \hline
  C\Cov(x_t, \hat{x}_t) & C\Var(x_t) C^{\top} + V
\end{array}
\bigg] \label{eqn:KF_grad2}\\
& G_t = \big[
\begin{array}{c|c}
  A \Cov(x_t, \hat{x}_t) & A\Var(x_t) C^{\top}
\end{array}\big] \label{eqn:KF_grad3}\\
& \Var(\hat{x}_{t+1}) = \pi_t \Psi_t \pi_t^{\top}, \ \Var(x_{t+1}) = A\Var(x_t)A^{\top} + W \nonumber\\
&\Cov(x_{t+1}, \hat{x}_{t+1}) = A\big[\Cov(x_t, \hat{x}_t) \mid  \Var(x_t)C^{\top}\big]\pi_t^{\top} \nonumber\\
&\Var(\hat{x}_0) = \Cov(x_0, \hat{x}_0) = \bar{x}_0\bar{x}_0^{\top}, \ \Var(x_0) = \bar{x}_0\bar{x}_0^{\top}+X_0. \nonumber
\end{align}
We next define the vanilla PG update as
\begin{align}\label{eqn:KF_update}
	\pi_h' = \pi_h - \eta_h\cdot \nabla_{\pi_h}\cJ_h(\pi_h). 
\end{align}
where $\eta_h > 0$ is a constant stepsize. When the exact PG in \eqref{eqn:PG_RHPG_KF} is not available, it can be estimated from samples of system trajectories using (two-point) zeroth-order optimization techniques as described in Algorithm \ref{alg:zeroth_order_GD}. Due to the landscape properties listed in Theorem \ref{theorem:quadratic}, global convergence and sample complexity of the PG update \eqref{eqn:KF_update} and its zeroth-order implementation naturally follow. We present them in the following propositions. 

\begin{algorithm}[t]
\caption{Two-Point Zeroth-Order Oracle}\label{alg:zeroth_order_GD}
\begin{algorithmic}[1]
\renewcommand{\algorithmicrequire}{\textbf{Input:}}
 \renewcommand{\algorithmicensure}{\textbf{Output:}}
 \REQUIRE Iteration index $h$, filter $\pi_h$, smoothing radius $r_h$
\STATE Uniformly sample $U_h$ from the surface of a unit sphere; Compute $\pi_h^{+} \leftarrow \pi_h + r_hU_h$ and $\pi_h^{-} \leftarrow \pi_h - r_hU_h$
\STATE Sample $x_0$ and compute $x_h$ and $\hat{x}_h$ by \eqref{eqn:inject1}-\eqref{eqn:inject2} using convergent filters $\{\pi_t\}_{t\in\{0, \cdots, h-1\}}$
\STATE Generate $x_{h+1}$ and $y_h$ using a simulator.
\STATE Compute estimates of $x_{h+1}$ using $\pi_h^{+}$ and $\pi_h^{-}$, denoted as $\hat{x}_{h+1}^{+}$ and $\hat{x}_{h+1}^{-}$, respectively 
\STATE Calculate objective values  $J_h(\pi_h^{+})$ and $J_h(\pi_h^{-})$ of   \eqref{eqn:induction} 

\STATE \textbf{return} $\tilde{\nabla}_{\pi_h}\cJ_h(\pi_h) \leftarrow \frac{n(m+n)}{2r_h}\cdot\big[J_h(\pi_h^{+})-J_h(\pi_h^{-})\big]\cdot U_h$
\end{algorithmic}
\end{algorithm}

\begin{proposition}\label{prop:global_conv}
	For all $h \in \{0, \cdots, N-1\}$ and a fixed $\epsilon > 0$, choose a constant stepsize $\eta_h \leq   1/\lambda_{\max}(H_h)$. Then, the PG update \eqref{eqn:KF_update} converges linearly to the unique minimum of \eqref{eqn:induction}. That is, such that $\|\pi_h^{T_h} - \tilde{\pi}_h^*\big\| \leq \epsilon$ after a total number of $T_h\sim\cO(\log(\epsilon^{-1}))$ iterations, where $\tilde{\pi}_h^*$ is the unique minimizing solution of  \eqref{eqn:induction} after setting filters for all $t \in \{0, \cdots, h-2\}$ to be those computed in the previous iterations and are $\epsilon$-close to the unique minimum, respectively. 
\end{proposition}

\begin{proposition}\label{prop:sample}
Choose the smoothing radius of the zeroth-order PG to satisfy $r_{h} \sim \cO(\epsilon)$ and the stepsize $\eta_{h}\sim \cO(\epsilon^2)$. Then, the zeroth-order PG update in Algorithm \ref{alg:zeroth_order_GD} converges after $T_h \sim \tilde{\cO}(\epsilon^{-2}\log(\frac{1}{\delta\epsilon^2}))$ iterations in the sense that $\|\pi_h^{T_h} - \tilde{\pi}_h^*\big\| \leq \epsilon$ with a probability of at least $1-\delta$.
\end{proposition}

Proposition \ref{prop:global_conv} is standard, and Proposition \ref{prop:sample} follows from the proof of Proposition 3.3 in \cite{zhang2023revisiting}. Combining Theorem \ref{theorem:KF_DP} with Proposition \ref{prop:sample}, we conclude that if we spend $\tilde{O}(\epsilon^{-2})$ samples in solving every one-step KF problem to $\cO(\epsilon)$-accuracy with a probability of $1-\delta$, for all $h \in \{0, \cdots, N-1\}$, then Algorithm \ref{alg:DPfilter} is guaranteed to output $\pi_{N-1}$ that is $\epsilon$-close to the infinite-horizon KF with a probability of at least $1-N\delta$.  The total sample complexity of the RHPG algorithm is thus $\tilde{O}(\epsilon^{-2})\cdot O(\log(\epsilon^{-1}))\sim \tilde{O}(\epsilon^{-2})$. This complexity result matches the complexity of applying RHPG to the LQR task \cite{zhang2023revisiting}.

\section{Numerical Experiment: Estimation of the Convection-Diffusion Model}\label{sec:sim}
We conducted numerical experiments to design state estimators for the one-dimensional convection-diffusion linear PDE\footnote{We open-source the code repository at \url{https://github.com/xiangyuan-zhang/LearningKF}}. The convection-diffusion equation models physical phenomena involving the transfer of particles, energy, or other quantities within a system due to convection and diffusion. These quantities are described by a continuous concentration function $c(\mathsf{x},t):\Omega \times \mathbb{R}^+ \rightarrow \mathbb{R}$, where $\mathsf{x}$ and $t$ represent spatial and temporal coordinates, respectively, and $\Omega \subset \mathbb{R}$ is the spatial domain of interest. The one-dimensional convection-diffusion equation can then be expressed as
\begin{align}
	\frac{\partial{c}}{\partial t} = \nu \frac{\partial^2 c}{\partial \mathsf{x}^2} - v \frac{\partial c}{\partial \mathsf{x}}, \label{eqn:cd_pde}
\end{align}
where $\nu$ is the diffusion coefficient and $v$ is the convection velocity; these scalar physical parameters characterize the strength of convection and diffusion, respectively. When $v = 0$, the convection-diffusion equation \eqref{eqn:cd_pde} reduces to the heat equation. As with any PDE, the convection-diffusion equation must be accompanied by initial and boundary conditions. Here, we consider the domain $\Omega = [0,1]$ with periodic boundary conditions, while the initial condition will be defined shortly.

The convection-diffusion equation can be solved numerically by discretizing space and time, resulting in a linear state-space model of the form \eqref{eqn:inf_LQE_dynamics}. To do so, we define a state vector $x_t \in \mathbb{R}^n$ that contains the values of $c$ at $n$ equally-spaced points in $\Omega$ and at time $t = k \Delta t$, where $n$ is even, $k = 0, 1, \dots$, and $\Delta t$ is the discrete time step. The dynamics governed by the convection-diffusion equation can then be approximated by a state-space model $x_{t+1} = A x_t$ with
\begin{align}\label{eqn:pde_A_matrix}
	A = \frac{1}{n} \cdot D^\dagger \diag (e^{-i v k_\mathsf{x} - \nu k_\mathsf{x}^2} \Delta t) D,
\end{align}
where $k_\mathsf{x} = 2 \pi [0, \dots, n/2-1, 0, -n/2+1, \dots, -1]$ is the vector of spatial wavenumbers, $i$ is the imaginary unit, $D$ is the discrete Fourier transform (DFT) matrix defined by $D_{pq} = e^{-2 \pi i (p-1)(q-1)/n}$, and its scaled conjugate transpose $D^\dagger/n$ is the inverse discrete Fourier transform (IDFT) matrix \cite{rao2018transform}. The matrix $A$ combines a spectral discretization of the spatial derivatives in the convection-diffusion equation, which takes advantage of the periodicity of the spatial domain with an exact temporal integration of its continuous-time dynamics. Such spectral evaluation of the derivatives enjoys exponential convergence properties \cite{trefethen1996finite}. As a result, even a small state dimension $n$ yields a state-space model that faithfully reproduces the dynamics of the convection-diffusion equation.

\begin{figure}[t]
\centering 
\includegraphics[width = 0.55\textwidth]{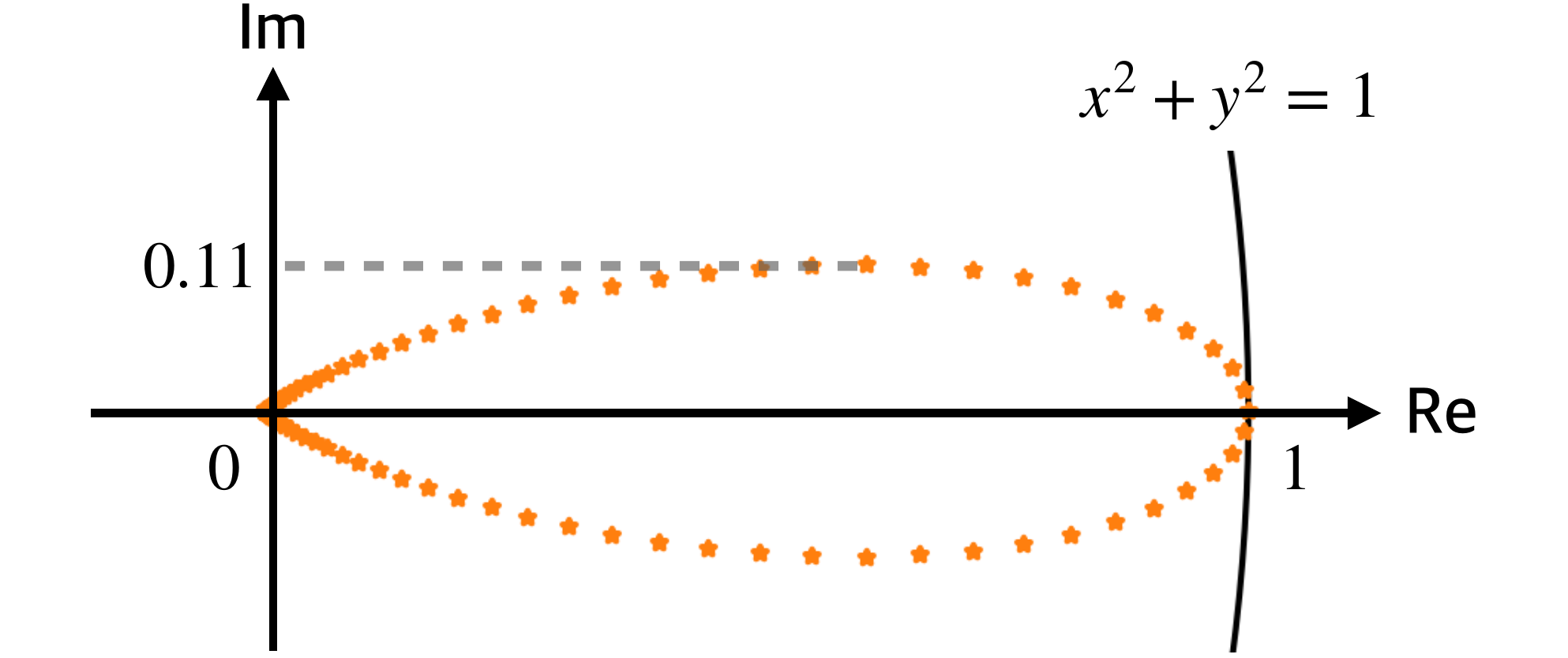}	
\caption{Eigen-spectrum of $A$ constructed in (\ref{eqn:pde_A_matrix}), where $\rho(A) = 1$ since the integrated concentration over the domain is conserved by the convection-diffusion equation.}\label{fig:eigen}
\end{figure}

To set up the numerical experiments, we chose $\Delta t = 0.05$, and set the dimensions of the state vector $x_t$ and measurement/observation vector to $n = 200$ and $m=5$, respectively. The five sensors were evenly distributed across the physical domain $\Omega$, where each sensor measures the (unscaled) values of the state at the corresponding location, subject to additive zero-mean Gaussian white noises. Moreover, we chose the diffusion coefficient to $\nu = 2\times 10^{-3}$, the convection velocity to $v = 5\times 10^{-2}$, the distribution of $x_0$ to be
\begin{align}\label{eqn:init_dist}
	x_0 \sim \cN\Big(&c(\mathsf{x},t=0), \frac{1}{16}\cdot\sin(2\pi \mathsf{x}) (\sin(2\pi \mathsf{x}))^{\top}\Big) \\
	&c(\mathsf{x},t=0) = \mathrm{sech}(10(\mathsf{x}-1/2)). \label{eqn:init_cond}
\end{align}
Furthermore, we set the covariance matrix of the measurement noise to $V = 10^{-1}\cdot\bI$, and the covariance matrix of the process to $W = 10^{-9}\cdot \bI$. Lastly, we set the covariance matrix of the additional noise $\theta_0$ in \eqref{eqn:inject1}-\eqref{eqn:inject2} to $\Theta = 10^{-2}\cdot \bI$.

\renewcommand\arraystretch{1.15}
\begin{algorithm}[t]
\caption{First-Order RHPG with \textsc{Adam}}\label{alg:analytic_gd}
\begin{algorithmic}[1]
\renewcommand{\algorithmicrequire}{\textbf{Input:}}
 \renewcommand{\algorithmicensure}{\textbf{Output:}}
 \REQUIRE Problem horizon $N$, $\Theta$
 \STATE Compute $\Delta$, $\Xi$ by \eqref{eqn:KF_grad1} 
\FOR{$h \in \{0, \cdots, N-1\}$}
	\STATE Compute $\Psi_h$, $G_h$ by \eqref{eqn:KF_grad2}-\eqref{eqn:KF_grad3}
	\STATE $i_h \leftarrow 0$, \ $m_h \leftarrow \bm{0}_{n \times (n+m)}$, $v_h \leftarrow \bm{0}_{(n+m) \times (n+m)}$ 
        \WHILE{$\|\nabla_{\pi_h}\| \geq 10^{-4}$}
        \STATE $\nabla_{\pi_h} \leftarrow 2 \big[\pi_h(\Psi_h + \Delta) - (G_h + \Xi)\big]$
         \STATE $\pi_h, m_h, v_h \leftarrow \textsc{Adam}(\pi_h, \nabla_{\pi_h}, m_h, v_h, 10^{-3}, i_h)$
         \STATE $i_h \leftarrow i_h + 1$
        \ENDWHILE
        	\IF{$h\neq N-1$}
        		\STATE $\pi_{h+1} \leftarrow \pi_h$  \textcolor{darkgray}{// Warm start the next iteration}
        	\ENDIF 
 \ENDFOR

\STATE \hspace{0.5cm}\textbf{return} $\pi_{N-1}$
\STATE 
\STATE {\textsc{Adam}}$(P, \nabla, m, v, \eta, i, \beta_1 = 0.9, \beta_2 = 0.999, \varsigma = 10^{-8})$
\STATE \hspace{0.5cm}$m \leftarrow \beta_1 m+ (1-\beta_1)\nabla$, ~$v \leftarrow \beta_2 v+ (1-\beta_2)(\nabla^{\top}\nabla)$
\STATE \hspace{0.5cm}$\hat{m} \leftarrow m/(1-\beta_1^i)$, \ $\hat{v} \leftarrow v/(1-\beta_2^i)$
\STATE \hspace{0.5cm}$P \leftarrow P - \eta\cdot \hat{m}(\hat{v}^{1/2} + \varsigma\bI)^{-1}$
\STATE \hspace{0.5cm}\textbf{return}  $P, m, v$
\end{algorithmic}
\end{algorithm}

After setting up the PDE environment, we applied the RHPG algorithm to learn the KF design with several different problem horizons ranging from $N=1$ to $N = 101$. We implemented the inner loop of the RHPG algorithm using the first-order vanilla PG update \eqref{eqn:KF_update} with the stepsize selected based on the \textsc{Adam} rule \cite{kingma2014adam}. We detail the procedure of the first-order RHPG algorithm in Algorithm \ref{alg:analytic_gd}. 

\begin{figure*}[t]
\centering 
\includegraphics[width = 0.87\textwidth]{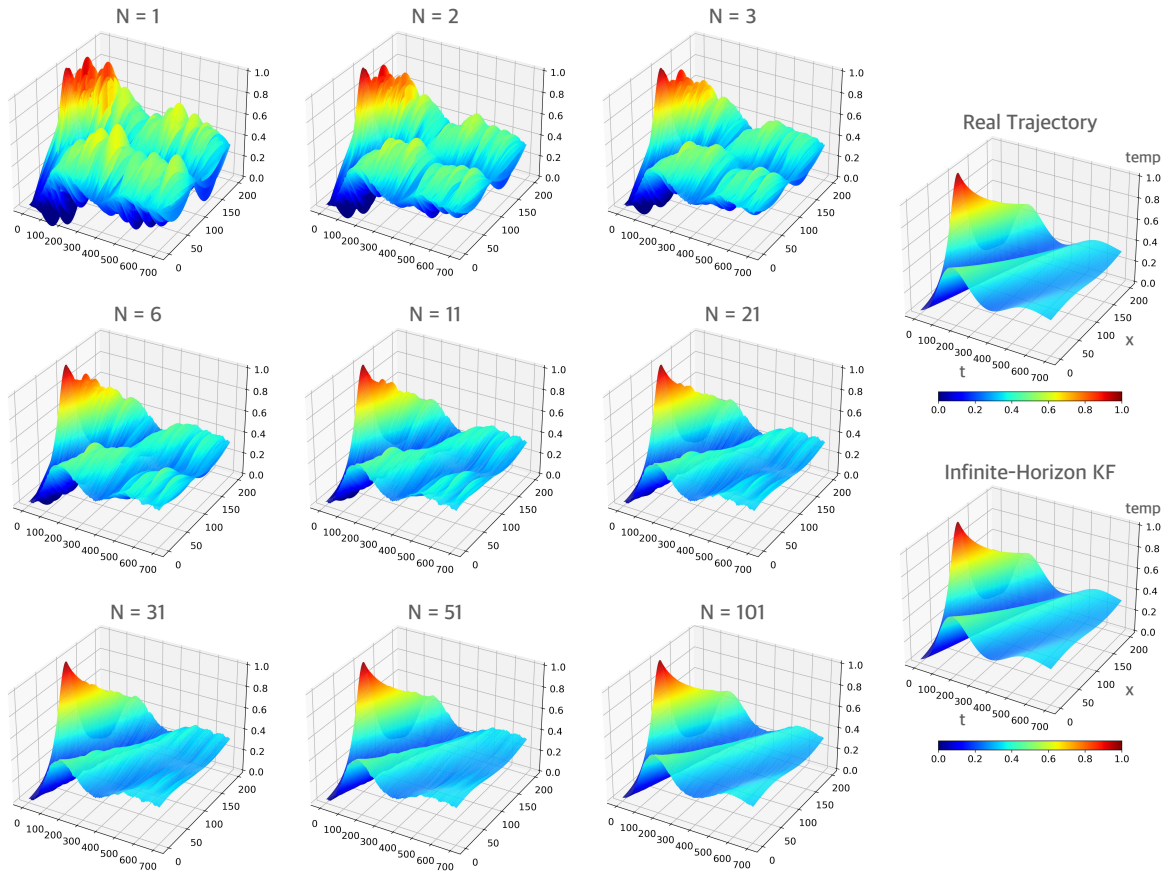}	
\caption{Visualizing the estimated trajectories by RL-based filters through RHPG and the model-based infinite-horizon KF (\ref{eqn:dynamic_filter})-(\ref{eqn:kalman_gain}), along with the ground truth system trajectory. The $x$, $y$, and $z$ axes of the three-dimensional plots represent time ($700$ total steps), space ($200$ dimensions), and the normalized concentration of the convection-diffusion equation. With large $N$, the learned filter converges globally to the infinite-horizon KF.}\label{fig:heat_rhpg}
\end{figure*}

In the first experiment, we generated a ground truth system trajectory with the length of $700$ discrete time steps and with a deterministic initial condition being $x_0 = c(\mathsf{x}, t=0)$ in \eqref{eqn:init_cond}. Then, we visualized the estimated state trajectories computed using the convergent filters from RHPG with several different input horizons $N \in \{1, 2, 3, 6, 11, 21, 31, 51, 101\}$), where for all filters the initial state estimate was set to $\hat{x}_0 = c(\mathsf{x},t=0)$. In Figure \ref{fig:heat_rhpg}, we provide a comparison of the estimated trajectories by the RL-based filters with the ground truth trajectory and the estimated trajectory generated by the infinite-horizon KF \eqref{eqn:dynamic_filter}-\eqref{eqn:kalman_gain}. One can observe from the first row of Figure \ref{fig:heat_rhpg} that when $N$ is small, the convergent filters of RHPG are myopic and can only predict system states over a short period since the convergences of RHPG in the initial iterations are toward the time-varying KFs in \eqref{eqn:kalman_gain_finite}. As we increase the number of iterations of RHPG, the filter adapts to new tasks with longer and longer problem horizons (cf., the illustration in Figure \ref{fig:landscape}). Lastly, one can observe from the third row of Figure \ref{fig:heat_rhpg} that when the problem horizon becomes sufficiently large, the RHPG algorithm converges globally to the infinite-horizon KF, which corroborates the theories developed in the paper.

\begin{figure}[t]
\centering 
\includegraphics[width = 0.55\textwidth]{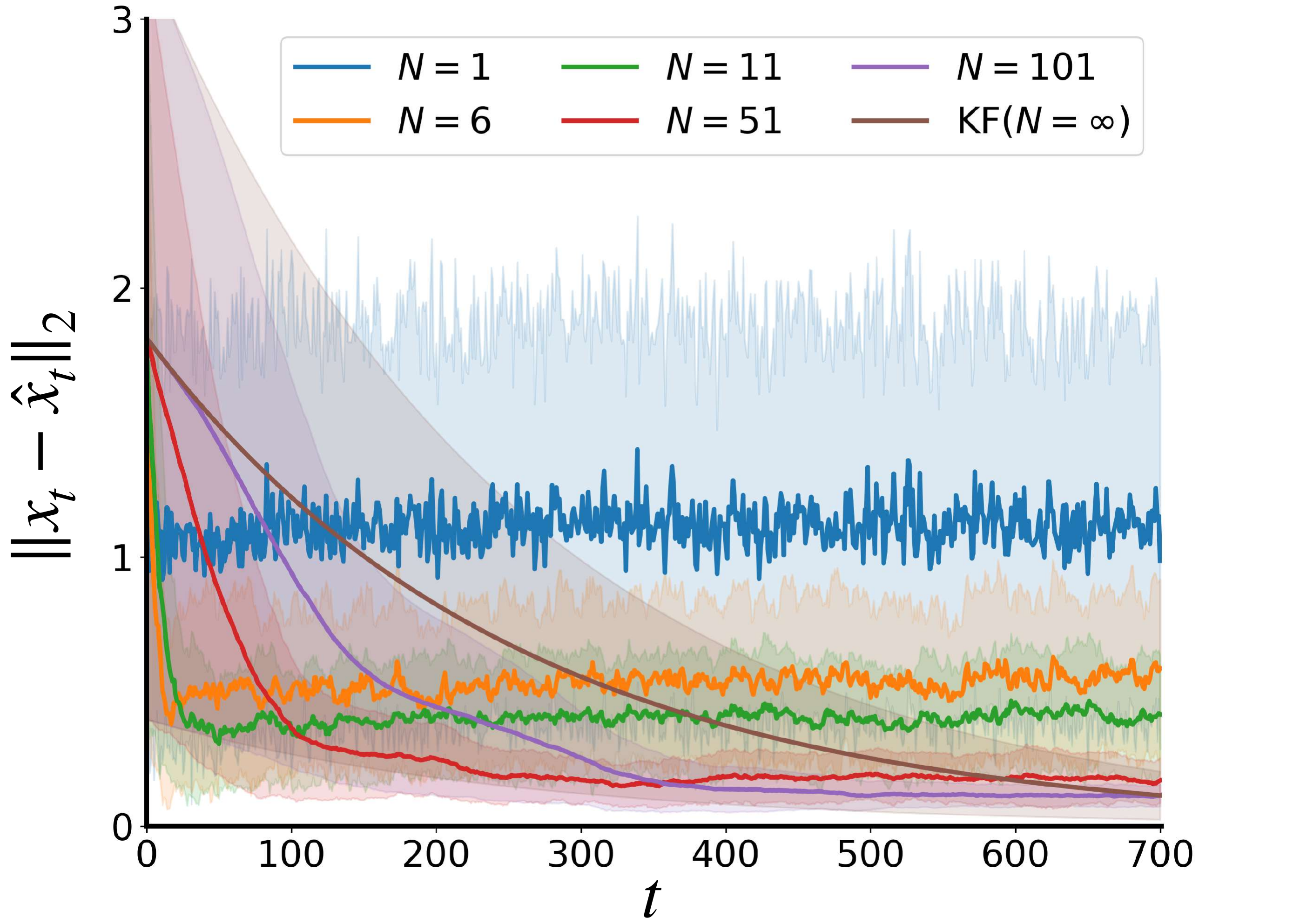}	
\caption{Average state estimation error $\|x_t-\hat{x}_t\|_2$ of RHPG filters and the model-based KF over $100$ different system trajectories. Each trajectory starts from a random initial condition and is affected by randomly sampled process and measurement noises.}\label{fig:kf_test}
\end{figure}

In the second experiment, we generated $100$ random system trajectories of the convection-diffusion equation for $700$ time steps, each starting from a randomly sampled initial condition according to \eqref{eqn:init_dist}. We then applied the convergent filters of RHPG as well as the model-based KF to estimate the $100$ system trajectories and plotted the average estimation error (i.e., $\|x_t-\hat{x}_t\|_2$) over time. As shown in Figure \ref{fig:kf_test}, the RHPG filters with a small problem horizon $N$ are myopic, resulting in lower estimation costs over the short time period. However, these myopic filters are not able to regulate the state estimation error in the asymptotic regime. When we run RHPG for a sufficiently large number of steps (e.g., $51$ and $101$), the learned filter performs well in the asymptotic regime due to the global convergence of RHPG to the infinite-horizon KF.

\section{Conclusion}\label{sec:conclude}
We have introduced the RHPG algorithm and provided rigorous analyses for its convergence and sample complexity in learning the infinite-horizon KF. RHPG is the first model-free PG algorithm with provable global convergence in learning the KF design, and it does not require any prior knowledge of the system for initialization and does not require the target system to be open-loop stable. We have validated our theories in extensive numerical experiments on a large-scale convection-diffusion model. On a higher level, our work has proposed developing RL algorithms specifically for control applications by utilizing classic control theory in the algorithmic design. The proposed paradigm enables certifying provable performance guarantees in a purely model-free setting, overcoming the nonconvex optimization landscape. Following this work and the dual theory to LQR \cite{zhang2023revisiting}, several ongoing and future research directions include designing and analyzing RHPG-type algorithms for LQG and $\cH_{\infty}$-robust filtering.
 
\section*{Acknowledgment}
X. Zhang and T. Ba\c{s}ar were supported in part by the US Army Research Laboratory (ARL) Cooperative Agreement W911NF-17-2-0181, and in part by the Army Research Office (ARO) MURI Grant AG285. S. Mowlavi and M. Benosman were supported solely by Mitsubishi Electric Research Laboratories. X. Zhang acknowledges helpful discussions with Bin Hu of UIUC in the early stage of the project and with Arvind Raghunathan of MERL. X. Zhang and T. Ba\c{s}ar acknowledge anonymous reviewers of ACC '23 for their helpful comments.

\small
\bibliographystyle{IEEEtran}
\bibliography{main}
\normalsize

\renewcommand{\theequation}{\thesection.\arabic{equation}}
\renewcommand{\thetheorem}{\thesection.\arabic{theorem}}

\appendix
\section{Proof of Theorem \ref{lemma:finite_approximation}}\label{proof:finite}
We first present a technical lemma due to \cite{de1989monotonicity}.
\begin{lemma} \label{lemma:comparison}
	Consider two RDEs 
\begin{align*}
	\Sigma_{t+1}^1 = W + A\Sigma_{t}^1A^{\top} - A\Sigma_{t}^1C^{\top}(V + C\Sigma_{t}^1C^{\top})^{-1}C\Sigma_{t}^1A^{\top}, \\
	\Sigma_{t+1}^2 = W + A\Sigma_{t}^2A^{\top} - A\Sigma_{t}^2C^{\top}(V + C\Sigma_{t}^2C^{\top})^{-1}C\Sigma_{t}^2A^{\top}.
\end{align*}
Then, the difference between the two solutions, $\tilde{\Sigma}_t:=\Sigma^2_t-\Sigma^1_t$, for all $t$, satisfies
\begin{align}
	\tilde{\Sigma}_{t+1} = \overline{A}_t\tilde{\Sigma}_t\overline{A}^{\top}_t - \overline{A}_t\tilde{\Sigma}_tC^{\top}(\tilde{V}_t + C\tilde{\Sigma}_tC^{\top})^{-1}C\tilde{\Sigma}_t\overline{A}^{\top}_t, 
\end{align}
where $\tilde{V}_t = V+C\Sigma^1_tC^{\top}$ and $\overline{A}_t:=A-A\Sigma^1_tC^{\top}\tilde{V}_t^{-1}C$.
\end{lemma}

Next, identify $\Sigma^1_t$ with $\Sigma^*$ and $\Sigma^2_t$ with $\Sigma^*_t$ in Lemma \ref{lemma:comparison}. Then, $\tilde{V}_t \hspace{-0.1em}=\hspace{-0.1em} V+C\Sigma^*C^{\top} \hspace{-0.1em}=:\hspace{-0.1em}\tilde{V}$ and $\overline{A}_t = A-A\Sigma^*C^{\top}\tilde{V}^{-1}C =: \overline{A}$ for all $t$. Invoking Lemma \ref{lemma:comparison} leads to
\begin{align}
	\tilde{\Sigma}_{t+1} &= \overline{A}\tilde{\Sigma}_t\overline{A}^{\top} - \overline{A}\tilde{\Sigma}_tC^{\top}(\tilde{V} + C\tilde{\Sigma}_tC^{\top})^{-1}C\tilde{\Sigma}_t\overline{A}^{\top} \label{eqn:RDE_conv_step1}\\
	&= \overline{A}\tilde{\Sigma}_t^{1/2}\big[\bI + \tilde{\Sigma}^{1/2}_tC^{\top}\tilde{V}^{-1}C\tilde{\Sigma}_t^{1/2}\big]^{-1}\tilde{\Sigma}_t^{1/2}\overline{A}^{\top} \nonumber \\
	&\leq \big[1+\lambda_{\min}(\tilde{\Sigma}^{1/2}_tC^{\top}\tilde{V}^{-1}C\tilde{\Sigma}_t^{1/2})\big]^{-1}\overline{A}\tilde{\Sigma}_t\overline{A}^{\top} =: \mu_t \overline{A}\tilde{\Sigma}_t\overline{A}^{\top},\label{eqn:RDE_conv_step2}
\end{align}
where $\tilde{\Sigma}_t^{1/2}$ denotes the unique psd square root of the psd matrix $\tilde{\Sigma}_t$, $0 < \mu_t \leq 1$ for all $t$, and $\overline{A}$ is the closed-loop matrix of the optimal infinite-horizon KF that has all its eigenvalue inside the unit circle (i.e., $\rho(\overline{A}) < 1$). Next, we use $\|\cdot\|_*$ to represent the $\Sigma^*$-induced matrix norm defined as $\|X\|^2_* := \max_{z \neq 0} \frac{z^{\top}X^{\top}\Sigma^*Xz}{z^{\top}\Sigma^*z}$. Then, we invoke Theorem 14.4.1 of \cite{hassibi1999indefinite}, where our $\tilde{\Sigma}_t$, $\overline{A}$ and $\Sigma^*$ correspond to $P_i - P^*$, $F_p$ and $W$ in \cite{hassibi1999indefinite}, respectively. By Theorem 14.4.1 of \cite{hassibi1999indefinite} and \eqref{eqn:RDE_conv_step2}, we obtain $\|\overline{A}\|_* < 1$ and given that $\mu_t \leq 1$, $\|\tilde{\Sigma}_{t+1}\|_* \leq \|\overline{A}\|^2_* \cdot \|\tilde{\Sigma}_t\|_*$. Therefore, the convergence rate is exponential in the sense that $\|\tilde{\Sigma}_t\|_* \leq \|\overline{A}\|_*^{2t}\cdot \|\tilde{\Sigma}_{0}\|_*$. Next, recall the condition number of a matrix $X$ is defined as $\kappa_X := \sigma_{\max}(X)/\sigma_{\min}(X)$. Then, the convergence of $\tilde{\Sigma}_t$ to the zero matrix in spectral norm can be characterized as
\begin{align*}
	\|\tilde{\Sigma}_t\| \leq \kappa_{\Sigma^*}\cdot \|\tilde{\Sigma}_t\|_* \leq \kappa_{\Sigma^*}\cdot\|\overline{A}\|_*^{2t}\cdot \|\tilde{\Sigma}_{0}\|_*.
\end{align*}
In other words, to ensure $\|\tilde{\Sigma}_N\| \leq \epsilon$, it suffices to require
\begin{align}\label{eqn:required_time}
	N \geq \frac{1}{2}\cdot \frac{\log\big(\frac{\|\tilde{\Sigma}_0\|_*\cdot \kappa_{\Sigma^*}}{\epsilon}\big)}{\log\big(\frac{1}{\|\overline{A}\|_*}\big)}.
\end{align}
Furthermore, since $(A, W^{1/2})$ is controllable and $(C, A)$ is observable, and if further $X_0 > \Sigma^*$, the closed-loop system at any time $t \geq 0$ is exponentially asymptotically stable such that the time-invariant (frozen) filter satisfies $\rho(A-A\Sigma^*_tC^{\top}(V+C\Sigma^*_tC^{\top})^{-1}C) < 1$ \cite{bitmead1985monotonicity, de1989monotonicity}.  Lastly, we show that the (monotonic) convergence of the filter gain to the time-invariant Kalman gain follows from the convergence of $\Sigma^*_t$ to $\Sigma^*$, which can be verified through: 
\begin{align}
	L^*_t - L^* &= A\Sigma^*_tC^{\top}\hspace{-0.1em}(V+C\Sigma^*_tC^{\top})^{\hspace{-0.1em}-1} \hspace{-0.2em}-\hspace{-0.15em} A\Sigma^*C^{\top}\hspace{-0.1em}(V+C\Sigma^*C^{\top})^{\hspace{-0.1em}-1} \nonumber\\
 	&=A\Sigma^* C^{\top}\big[ (V + C\Sigma^*_t C^{\top})^{-1} - (V + C\Sigma^* C^{\top})^{-1}\big] + A(\Sigma^*_t - \Sigma^*) C^{\top}(V + C\Sigma^*_t C^{\top})^{-1}\nonumber\\
 	&= A\Sigma^* C^{\top}(V + C\Sigma^* C^{\top})^{-1}C(\Sigma^*-\Sigma^*_t)C^{\top}(V + C\Sigma^*_t C^{\top})^{-1} -A(\Sigma^*- \Sigma^*_t) C^{\top}(V + C\Sigma^*_t C^{\top})^{-1} \nonumber\\
 	&=(L^*C - A)(\Sigma^* - \Sigma^*_t) C^{\top}(V + C\Sigma^*_tC^{\top})^{-1}. \label{eqn:ldiff}
\end{align}
Hence, we have $\|L^*_t - L^*\| \leq \frac{\|\overline{A}\|\cdot \|C\|}{\lambda_{\min}(V)}\cdot \|\Sigma^*_t - \Sigma^*\|$. Substituting $\epsilon$ in \eqref{eqn:required_time} with $\frac{\epsilon\cdot\lambda_{\min}(V)}{\|\overline{A}\|\cdot\|C\|}$ and identifying that $\overline{A}$ is exactly $A_L^*$ completes the proof.

\section{Proof of Theorem \ref{theorem:KF_DP}}\label{proof:KF_DP}
To prove $\big\|[\tilde{A}_{L_{N-1}} \ \tilde{B}_{L_{N-1}}] - [A^*_L \ B^*_L]\big\| \leq \epsilon$, it suffices to bound the error between the approximated filter gain $\tilde{L}_{N-1}$ and the exact Kalman gain $L^*$ as in \eqref{eqn:kalman_gain}. First, according to Theorem \ref{lemma:finite_approximation}, we select
\begin{align}\label{eqn:N_choice}
	N = \frac{1}{2}\cdot \frac{\log\big(\frac{2\|X_0-\Sigma^*\|_*\cdot\kappa_{\Sigma^*}\cdot \|A_L^*\|\cdot\|C\|} {\epsilon\cdot\lambda_{\min}(V)}\big)}{\log\big(\frac{1}{\|A_L^*\|_*}\big)} + 1.
\end{align} 
which ensures that $L^*_{N-1}$ is stabilizing and $\|L^*_{N-1} - L^*\| \leq \epsilon/2$. Then, it remains to show that Algorithm \ref{alg:DPfilter} returns a filter $\tilde{L}_{N-1}$ such that $\|\tilde{L}_{N-1} - L^*_{N-1}\| \leq \epsilon/2$. 

Recall that the FRDE is the following forward iteration starting with $\Sigma^*_0 = X_0 > 0$:
\begin{align}
	\Sigma^*_{t+1} &= A\Sigma^*_{t}A^{\hspace{-0.1em}\top} \hspace{-0.2em}-\hspace{-0.1em} A\Sigma^*_t C^{\top}\hspace{-0.1em}(V \hspace{-0.1em}+\hspace{-0.1em} C\Sigma^*_t C^{\top})^{\hspace{-0.1em}-1}C\Sigma^*_t A^{\top} \hspace{-0.1em}+\hspace{-0.1em} W \label{eqn:standard_RDE_appen}\\
	&= (A-L^*_tC)\Sigma^*_tA^{\top} + W \label{eqn:filter_riccati_finite_appen}\\
	&=(A-L^*_tC)\Sigma^*_t(A-L^*_tC)^{\top} + L^*_tV(L^*_t)^{\top}+W \label{eqn:filter_RDE_Lya}. 
\end{align}
 Moreover, for an arbitrary $L_t$, it holds that:
 \begin{align}
 	\Sigma_{t+1} = (A-L_tC)\Sigma_t(A-L_tC)^{\top} + L_tVL_t^{\top}+W. \label{eqn:filter_lyapunov}
 \end{align}
Furthermore, for clarity of the proof, we define/recall:
\begin{align*}
 	&L^*_t\text{: Exact Kalman gain at time $t$ defined in \eqref{eqn:kalman_gain_finite}}\\
 	&\tilde{L}_t^*\text{: Optimal gain of the current cost-to-come function,}\\
 	&\hspace{1.6em} \text{absorbing errors in prior steps} \\
 	&\tilde{L}_t\text{: An approximation of $\tilde{L}_t^*$ obtained by applying \eqref{eqn:KF_update}}\\
 	&\delta_t:=\tilde{L}_t-\tilde{L}_t^* \text{: Policy optimization error at time $t$} \\
 	&\tilde{\Sigma}^*_{t+1} \text{: Solution generated by \eqref{eqn:filter_RDE_Lya} with $L^*_t = \tilde{L}^*_t$ and $\Sigma^*_t = \tilde{\Sigma}_t$.}
 \end{align*}

We argue that $\|\tilde{L}_{N-1} - L^*_{N-1}\| \leq \epsilon/2$ can be achieved by carefully controlling $\delta_t$ for all $t$. At $t=N-1$, it holds that
\begin{align*}
	\|\tilde{L}_{N-1} - L^*_{N-1}\| &\leq \|\tilde{L}^*_{N-1} - L^*_{N-1}\| + \|\delta_{N-1}\|,
\end{align*}
where from \eqref{eqn:ldiff} we have
\begin{align*}
	\tilde{L}^*_{N-1} - L^*_{N-1} = (L^*_{N-1}C - A)(\Sigma^*_{N-1} - \tilde{\Sigma}_{N-1})C^{\top} \cdot(V +C\tilde{\Sigma}_{N-1}C^{\top})^{-1}.
\end{align*}
As a result, it holds that
\begin{align}
	\hspace{-0.3em}\|\tilde{L}^*_{N-1} \hspace{-0.1em}-\hspace{-0.1em} L^*_{N-1}\| \leq \frac{\|A_{L_{N-1}}^*\|\|C\|}{\lambda_{\min}(V)}\cdot\|\Sigma^*_{N-1} \hspace{-0.1em}-\hspace{-0.1em} \tilde{\Sigma}_{N-1}\| \label{eqn:laststep_req}.
\end{align}
Define the helper constant
\begin{align*}
	C_1:= \frac{\varphi\|C\|}{\lambda_{\min}(V)} > 0, \quad \varphi := \max_{t\in\{0, \cdots, N-1\}}\|A^*_{L_t}\|.
\end{align*}
Next, require $\|\delta_{N-1}\| \leq \frac{\epsilon}{4}$ and $\|\tilde{L}^*_{N-1} - L^*_{N-1}\| \leq \frac{\epsilon}{4}$ to fulfill $\|\tilde{L}_{N-1} - L^*_{N-1}\| \leq \frac{\epsilon}{2}$. By \eqref{eqn:laststep_req}, this is equivalent to requiring
\begin{align}\label{eqn:sigmareq1}
	\|\Sigma^*_{N-1} - \tilde{\Sigma}_{N-1}\| \leq \frac{\epsilon}{4 C_1}.
\end{align}
Subsequently, by \eqref{eqn:filter_lyapunov}, we have
\begin{align}\label{eqn:sigmalast}
	\Sigma^*_{N-1} \hspace{-0.1em}-\hspace{-0.1em} \tilde{\Sigma}_{N-1} = (\Sigma^*_{N-1} \hspace{-0.1em}-\hspace{-0.1em} \tilde{\Sigma}^*_{N-1}) \hspace{-0.1em}+\hspace{-0.1em} (\tilde{\Sigma}^*_{N-1} \hspace{-0.1em}-\hspace{-0.1em} \tilde{\Sigma}_{N-1}),
\end{align}
where the first difference term on the RHS of \eqref{eqn:sigmalast} is
\begin{align}
&\Sigma^*_{N-1} - \tilde{\Sigma}^*_{N-1} = (A-L^*_{N-2}C)\Sigma^*_{N-2}A^{\top} - (A-\tilde{L}^*_{N-2}C)\tilde{\Sigma}_{N-2}A^{\top} \nonumber\\
&= (A-L^*_{N-2}C)(\Sigma^*_{N-2} -\tilde{\Sigma}_{N-2})A^{\top} + (\tilde{L}^*_{N-2} - L^*_{N-2})C\tilde{\Sigma}_{N-2}A^{\top}.\label{eqn:sigmalast2}
\end{align}
Moreover, the second term on the RHS of \eqref{eqn:sigmalast} is 
\small
\begin{align}
	&\hspace{1em}\tilde{\Sigma}^*_{N-1} - \tilde{\Sigma}_{N-1} \nonumber\\
	&\hspace{-0.3em}= (A-\tilde{L}^*_{N-2}C)\tilde{\Sigma}_{N-2}(A-\tilde{L}^*_{N-2}C)^{\top} + \tilde{L}^*_{N-2}V(\tilde{L}^*_{N-2})^{\top}  - (A-\tilde{L}_{N-2}C)\tilde{\Sigma}_{N-2}(A-\tilde{L}_{N-2}C)^{\top} - \tilde{L}_{N-2}V(\tilde{L}_{N-2})^{\top} \nonumber\\
	&\hspace{-0.3em}= A\tilde{\Sigma}_{N-2}C^{\top}\tilde{L}_{N-2}^{\top} \hspace{-0.1em}-\hspace{-0.1em} A\tilde{\Sigma}_{N-2}C^{\top}(\tilde{L}^*_{N-2})^{\top} \hspace{-0.2em}+\hspace{-0.2em} \tilde{L}_{N-2}C\tilde{\Sigma}_{N-2}A^{\top} \nonumber\\
	&\hspace{1em}+ \tilde{L}^*_{N-2}(C\tilde{\Sigma}_{N-2}C^{\top} + V)(\tilde{L}^*_{N-2})^{\top} - \tilde{L}_{N-2}(C\tilde{\Sigma}_{N-2}C^{\top} + V)\tilde{L}_{N-2}^{\top} - \tilde{L}^*_{N-2}C\tilde{\Sigma}_{N-2}A^{\top} \nonumber\\
	&\hspace{-0.3em}=\tilde{L}^*_{N-2}(C\tilde{\Sigma}_{N-2}C^{\top} + V)(\tilde{L}^*_{N-2})^{\top}- A\tilde{\Sigma}_{N-2}C^{\top}(\tilde{L}^*_{N-2})^{\top} - \tilde{L}^*_{N-2}C\tilde{\Sigma}_{N-2}A^{\hspace{-0.1em}\top}  \hspace{-0.2em}+\hspace{-0.2em} A\tilde{\Sigma}_{N-2}C^{\hspace{-0.1em}\top}\hspace{-0.2em}(C\tilde{\Sigma}_{N-2}C^{\hspace{-0.1em}\top}\hspace{-0.2em}+\hspace{-0.2em}V)^{\hspace{-0.1em}-1}C\tilde{\Sigma}_{N-2}A^{\hspace{-0.1em}\top}\nonumber\\
	& -(\tilde{L}_{N-2} - A\tilde{\Sigma}_{N-2}C^{\top}(C\tilde{\Sigma}_{N-2}C^{\top}\hspace{-0.1em}+\hspace{-0.1em}V)^{-1})(C\tilde{\Sigma}_{N-2}C^{\top}\hspace{-0.1em}+\hspace{-0.1em}V)\cdot(\tilde{L}_{N-2} - A\tilde{\Sigma}_{N-2}C^{\top}(C\tilde{\Sigma}_{N-2}C^{\top}\hspace{-0.1em}+\hspace{-0.1em}V)^{-1})^{\top} \label{eqn:sigmadiff1}\\
	&\hspace{-0.3em}= (A\tilde{\Sigma}_{N-2}C^{\top}(C\tilde{\Sigma}_{N-2}C^{\top}\hspace{-0.1em}+\hspace{-0.1em}V)^{-1} \hspace{-0.1em}-\hspace{-0.1em} \tilde{L}^*_{N-2})(C\tilde{\Sigma}_{N-2}C^{\top}\hspace{-0.2em}+\hspace{-0.2em}V)\cdot(A\tilde{\Sigma}_{N-2}C^{\top}(C\tilde{\Sigma}_{N-2}C^{\top}+V)^{-1}-\tilde{L}^*_{N-2})^{\top}\nonumber\\
	&-(\tilde{L}_{N-2} - A\tilde{\Sigma}_{N-2}C^{\top}(C\tilde{\Sigma}_{N-2}C^{\top}\hspace{-0.2em}+\hspace{-0.2em}V)^{-1})(C\tilde{\Sigma}_{N-2}C^{\top}\hspace{-0.2em}+\hspace{-0.2em}V)\cdot(\tilde{L}_{N-2} - A\tilde{\Sigma}_{N-2}C^{\top}(C\tilde{\Sigma}_{N-2}C^{\top}\hspace{-0.1em}+\hspace{-0.1em}V)^{-1})^{\top}, \label{eqn:sigmadiff2}
\end{align}
\normalsize
where \eqref{eqn:sigmadiff1} is due to completion of squares. Substituting $\tilde{L}^*_{N-2} = A\tilde{\Sigma}_{N-2}C^{\top}(C\tilde{\Sigma}_{N-2}C^{\top}+V)^{-1}$ into \eqref{eqn:sigmadiff2} leads to
\begin{align}\label{eqn:sigmadiff3}
	\tilde{\Sigma}^*_{N-1} - \tilde{\Sigma}_{N-1} = -\delta_{N-2}(C\tilde{\Sigma}_{N-2}C^{\top}+V)\delta_{N-2}^{\top}.
\end{align}
Thus, combining \eqref{eqn:sigmalast}, \eqref{eqn:sigmalast2}, and \eqref{eqn:sigmadiff3} yields
\begin{align}
	&\hspace{1em}\|\Sigma^*_{N-1} - \tilde{\Sigma}_{N-1}\| \nonumber\\
	&\hspace{-0.3em}\leq \|\Sigma^*_{N-2} -\tilde{\Sigma}_{N-2}\|\varphi\|A\|  + \|\tilde{L}^*_{N-2} - L^*_{N-2}\|\|C\|\|\tilde{\Sigma}_{N-2}\|\|A\| + \|\delta_{N-2}\|^2\|C\tilde{\Sigma}_{N-2}C^{\top}+V\| \nonumber \\
	&\hspace{-0.3em}\leq\|A\|\cdot[\varphi+C_1\cdot\|C\|\cdot\|\tilde{\Sigma}_{N-2}\|]\cdot\|\Sigma^*_{N-2} - \tilde{\Sigma}_{N-2}\| + \|\delta_{N-2}\|^2\|C\tilde{\Sigma}_{N-2}C^{\top}+V\| \label{eqn:sigma_contraction},
\end{align}
where the last inequality follows from \eqref{eqn:laststep_req}. Now, require
\begin{align}
	\|\Sigma^*_{N-2} - \tilde{\Sigma}_{N-2}\| &\leq \frac{\epsilon}{4 C_1C_2} , \quad \|\delta_{N-2}\| \leq  \frac{1}{2}\sqrt{\frac{\epsilon}{C_1C_3}}, \label{eqn:sigmareq2}
\end{align}
where $C_2$ and $C_3$ are positive constants defined as
\begin{align*}
	C_2 := 2\|A\|\cdot\big[\varphi+C_1\cdot\|C\|\cdot\big(\|X_0\| + \|\Sigma^*\|\big)\big] >0, \quad C_3 := 2\big[\|V\| + \|C\|^2\big(\|X_0\|+ \|\Sigma^*\|\big)\big] > 0.
\end{align*}
Then, condition \eqref{eqn:sigmareq2} is sufficient for \eqref{eqn:sigmareq1} (and thus for $\|\tilde{L}_{N-1} - L^*_{N-1}\| \leq \epsilon/2$) to hold. Subsequently, we can propagate the requirements in \eqref{eqn:sigmareq2} backward in time. Specifically, we iteratively apply the arguments in \eqref{eqn:sigma_contraction} (i.e., by plugging quantities with subscript $t$ into the LHS of \eqref{eqn:sigma_contraction} and plugging quantities with subscript $t-1$ into the RHS of \eqref{eqn:sigma_contraction}) to obtain  the result that if at all $t \in \{1, \cdots, N-2\}$, we require
\begin{align}
	\hspace{-0.5em}\|\Sigma^*_{t} - \tilde{\Sigma}_{t}\| \leq \frac{\epsilon}{4 C_1C_2^{N-t-1}}, \ \|\delta_{t}\| \leq \frac{1}{2}\sqrt{\frac{\epsilon}{C_1C_2^{N-t-2}C_3}}\label{eqn:sigmareq_allt},
\end{align}
then \eqref{eqn:sigmareq2} holds true and therefore \eqref{eqn:sigmareq1} is satisfied. 

We now compute the required accuracy for $\delta_0$. As illustrated in Figure \ref{fig:proof_sketch}, we have $\Sigma^*_1 = \tilde{\Sigma}^*_1$ because 
\begin{align*}
	\tilde{\Sigma}^*_1 &= (A-\tilde{L}^*_0C)\Sigma^*_0(A-\tilde{L}^*_0C)^{\top} + \tilde{L}^*_0V(\tilde{L}^*_0)^{\top}+W = (A-L^*_0C)\Sigma^*_0(A-L^*_0C)^{\top} + L^*_0V(L^*_0)^{\top}+W = \Sigma^*_1,
\end{align*}
where the second equality is due to $\tilde{L}^*_0 = L^*_0$ since there are no prior computational errors yet at $t=0$. By \eqref{eqn:sigma_contraction}, the distance between $\Sigma^*_1$ and $\tilde{\Sigma}_1$ can be bounded as 
\begin{align*}
	\|\Sigma^*_1 - \tilde{\Sigma}_1\| = \|\tilde{\Sigma}^*_1 - \tilde{\Sigma}_1\| \leq \|\delta_{0}\|^2\cdot C_3.
\end{align*}
To fulfill the requirement \eqref{eqn:sigmareq_allt} for $t=1$, which is $\|\Sigma^*_{1} - \tilde{\Sigma}_{1}\| \leq \frac{\epsilon}{4 C_1C_2^{N-2}}$, it suffices to let
\begin{align}
	\|\delta_0\| \leq \frac{1}{2}\sqrt{\frac{\epsilon}{C_1C_2^{N-2}C_3}}. \label{eqn:delta0_req}
\end{align}

Lastly, we analyze the worst-case complexity of the proposed algorithm by computing, at the most stringent case, the required size of $\|\delta_t\|$. When $C_2 \leq 1$, the most stringent dependence of $\|\delta_t\|$ on $\epsilon$ happens at $t=N-1$, which is of the order $\cO(\epsilon)$, and the dependences on system parameters (through the dependence on constants $C_1, C_2$ and $C_3$) are polynomial. We argue that if $C_2 > 1$, then the requirement on $\|\delta_{N-1}\|$ is still the most stringent one. This is because $\|\delta_0\| \hspace{-0.1em}\leq\hspace{-0.1em} \|\delta_t\|$ for all $t \hspace{-0.1em}\in\hspace{-0.1em} \{1, \cdots, N\hspace{-0.1em}-\hspace{-0.1em}2\}$ and by \eqref{eqn:delta0_req}, we have
\begin{align}\label{eqn:delta0}
	\|\delta_0\| \sim \cO\Big(\sqrt{\frac{\epsilon}{C_2^{N-2}}}\Big).
\end{align}
Since we require $N$ to satisfy \eqref{eqn:N_choice}, the dependence of $\|\delta_0\|$ on $\epsilon$ in \eqref{eqn:delta0} becomes $\|\delta_0\| \sim \cO(\epsilon^{\frac{3}{4}})$, which is milder than that of $\|\delta_{N-1}\|$. Therefore, it suffices to require the most stringent error bound for all $t$, which is $\|\delta_t\| \sim \cO(\epsilon)$, to reach the $\epsilon$-neighborhood of the infinite-horizon KF. Lastly, for $\tilde{A}_{L_{N-1}}$ to be stable, it suffices to let $\epsilon$ to be sufficiently small such that $\epsilon < 1 - \|A_L^*\|_* \Longrightarrow \|\tilde{A}_L\|_* < 1$. This completes the proof. 

\section{Proof of Theorem \ref{theorem:quadratic}}\label{proof:quadratic}
We first introduce a standard result in linear algebra, see Eq. (7.7.5) and Theorem 7.7.7 on pages 495-496 of \cite{horn2012matrix}.
\begin{lemma}\label{lemma:pdness}
	For any symmetric matrix
	\begin{align*}
		X = \begin{bmatrix}
			X_{11} & X_{12} \\ X_{12}^{\top} & X_{22}
		\end{bmatrix} \in \RR^{(n+m)\times(n+m)},
	\end{align*}
	$X$ is pd if and only if $X_{11} > 0$ and $X_{22} - X_{12}^{\top}(X_{11})^{-1}X_{12} > 0$, where $X_{11} \in \RR^{n\times n}$, $X_{22} \in \RR^{m\times m}$, and $X_{12} \in \RR^{n\times m}$.
\end{lemma}

We start by analyzing the first iteration of the RHPG algorithm. The objective function for the one-step KF problem is defined as
\begin{align*}
	\cJ_0:= \EE_{x_0, w_0, v_0} \Big\{ \sum_{t=0}^{1}(x_t-\hat{x}_t)^{\top}(x_t-\hat{x}_t)\Big\}.
\end{align*}
The one-step static estimation problem, when formulated as a policy optimization problem, can be represented as 
\begin{align*}
	(\textbf{P1}) \ \min_{A_{L_0}, B_{L_0}} \ \cJ_0, \quad \text{ s.t. } \ x_1 = Ax_0 + w_0, \quad \hat{x}_1 = A_{L_0}\hat{x}_0 + B_{L_0}y_0, \quad x_0 = \cN(\bar{x}_0, X_0), \quad \hat{x}_0 = \bar{x}_0
\end{align*}
The Hessian matrix for the quadratic program $(\textbf{P1})$ is 
\begin{align*}
	H_{(\textbf{P1})} = \begin{bmatrix}
		\frac{\partial^2 \cJ_0}{A_{L_0}^2} & \frac{\partial^2 \cJ_0}{A_{L_0}\cdot B_{L_0}} \vspace{0.2em}\\ \frac{\partial^2 \cJ_0}{B_{L_0} \cdot A_{L_0}} & \frac{\partial^2 \cJ_0}{B_{L_0}^2}
	\end{bmatrix} = \begin{bmatrix}
		\bar{x}_0\bar{x}_0^{\top} & \bar{x}_0\bar{x}_0^{\top}C^{\top} \\ C\bar{x}_0\bar{x}_0^{\top} & C(\bar{x}_0\bar{x}_0^{\top}\hspace{-0.15em}+\hspace{-0.15em}X_0)C^{\top} \hspace{-0.25em}+\hspace{-0.15em} V
	\end{bmatrix}.
\end{align*}
Since $\bar{x}_0\bar{x}_0^{\top}$ is not pd (only one eigenvalue is positive) and according to Lemma \ref{lemma:pdness}, the quadratic program $(\textbf{P1})$ is not strictly convex, there exist multiple (in fact, infinite number of) stationary points that are equally good, and all of these stationary point policies attain the minimum value of (\textbf{P1}). However, we know that the one-step filter, denoted as $(A-L_0^*C, L_0^*)$, consists of one of such stationary points. 

To convexify (\textbf{P1}), we sample $\theta_0 \sim \cN(\bm{0}, \Theta)$ and set $\hat{x}_0 = \bar{x}_0 + \theta_0$, where $\theta_0$ is independent to $x_0$, $w_0$, and $v_0$. Concretely, we define the second quadratic program (\textbf{P2}) as follows:
\begin{align*}
	&(\textbf{P2}) \ \min_{A_{L_0}, B_{L_0}} \ \EE_{x_0, w_0, v_0, \theta_0} \Big\{ \sum_{t=0}^{1}(x_t-\hat{x}_t)^{\top}(x_t-\hat{x}_t)\Big\} \\ 
	\text{ s.t. } \  &x_1 = Ax_0 + w_0, \quad \hat{x}_1 = A_{L_0}\hat{x}_0 + B_{L_0}y_0, \quad x_0 = \cN(\bar{x}_0, X_0), \quad \hat{x}_0 = \bar{x}_0 + \theta_0
\end{align*}
In other words, the initial state estimate now satisfies $\hat{x}_0 \sim \cN(\bar{x}_0, \Theta)$. Then, (\textbf{P2}) can be shown to be equivalent to 
\begin{align}
	&\min_{A_{L_0}, B_{L_0}} \ \EE_{x_0, \theta_0}\Big\{ (x_0-\bar{x}_0 - \theta_0)^{\top}(x_0 -\bar{x}_0 - \theta_0)\Big\}\nonumber\\
	&\hspace{1em}+\EE_{x_0, \theta_0, w_0, v_0} \Big\{ (x_1-A_{L_0}(\bar{x}_0+\theta_0)-B_{L_0}y_0)^{\top}(x_1-A_{L_0}(\bar{x}_0+\theta_0)-B_{L_0}y_0)\Big\}\nonumber\\
	&\equiv \min_{A_{L_0}, B_{L_0}}\ \cJ_0 + \EE_{\theta_0}\Big\{\theta_0^{\top}(\bI + A_{L_0}^{\top}A_{L_0})\theta_0\Big\} \equiv \min_{A_{L_0}, B_{L_0}}\ \cJ_0 + \Tr\Big((\bI + A_{L_0}^{\top}A_{L_0})\Theta\Big). \label{eqn:regularizer}
\end{align}
Notice that $\Tr((\bI + A_{L_0}^{\top}A_{L_0})\Theta)$ in \eqref{eqn:regularizer} can be identified as an additional regularization term. The Hessian matrix becomes
\begin{align*}
		H_{(\textbf{P2})} = \begin{bmatrix}
		\bar{x}_0\bar{x}_0^{\top} + \Theta & \bar{x}_0\bar{x}_0^{\top}C^{\top} \\ C\bar{x}_0\bar{x}_0^{\top} & C(\bar{x}_0\bar{x}_0^{\top}\hspace{-0.15em}+\hspace{-0.15em}X_0)C^{\top} \hspace{-0.25em}+\hspace{-0.15em} V
	\end{bmatrix}.
\end{align*}
Since $\bar{x}_0\bar{x}_0^{\top} + \Theta > 0$ and 
\begin{align*}
	&C(\bar{x}_0\bar{x}_0^{\top}+X_0)C^{\top} + V - C\bar{x}_0\bar{x}_0^{\top}(\bar{x}_0\bar{x}_0^{\top} + \Theta)^{-1}\bar{x}_0\bar{x}_0^{\top}C^{\top}\\
	& \geq C(\bar{x}_0\bar{x}_0^{\top}+X_0)C^{\top} + V - C\bar{x}_0\bar{x}_0^{\top}(\bar{x}_0\bar{x}_0^{\top})^{\dagger}\bar{x}_0\bar{x}_0^{\top}C^{\top} = CX_0C^{\top} + V > 0,
\end{align*}
where we have used $\dagger$ to denote the Moore–Penrose matrix inverse, by Lemma \ref{lemma:pdness}, adding $\theta_0$ to the initial estimate $\hat{x}_0$ helps ensuring that (\textbf{P2}) is strictly convex. However, the solution to (\textbf{P2}) will not be identical to the solution to (\textbf{P1}) due to the extra regularization term. 

To achieve unbiased convexification, we inject $\theta_0$ to both $\hat{x}_0$ and $x_0$, which results in the following optimization problem
\begin{align*}
	&(\textbf{P3}) \ \min_{A_{L_0}, B_{L_0}} \ \EE_{x_0, w_0, v_0, \theta_0} \Big\{ \sum_{t=0}^{1}(x_t-\hat{x}_t)^{\top}(x_t-\hat{x}_t)\Big\} \\ 
	\text{ s.t. } \  &x_1 = Ax_0 + w_0, \quad \hat{x}_1 = A_{L_0}\hat{x}_0 + B_{L_0}y_0, \quad x_0 = \cN(\bar{x}_0, X_0) + \theta_0, \quad \hat{x}_0 = \bar{x}_0 + \theta_0
\end{align*}
Let $\Lambda_0 := A-B_{L_0}C-A_{L_0}$. Then, we can show that (\textbf{P3}) is equivalent to
\begin{align}
 &\min_{A_{L_0}, B_{L_0}} \EE_{x_0, \theta_0}\Big\{\hspace{-0.15em} (x_0-\bar{x}_0 +\theta_0 - \theta_0)^{\hspace{-0.15em}\top}\hspace{-0.15em}(x_0-\bar{x}_0+\theta_0 - \theta_0)\hspace{-0.15em}\Big\}\nonumber\\
	&\hspace{1.5em}+\EE \Big\{(Ax_0 + w_0-A_{L_0}\bar{x}_0 -B_{L_0}(Cx_0+v_0)+ \Lambda_0\theta_0)^{\top}(Ax_0 + w_0-A_{L_0}\bar{x}_0 -B_{L_0}(Cx_0+v_0)+ \Lambda_0\theta_0\Big\}\nonumber\\
	&\equiv \min_{A_{L_0}, B_{L_0}} \cJ_0 +\EE_{\theta_0}\Big\{\theta_0^{\top}\Lambda_0^{\top}\Lambda_0\theta_0\Big\}\equiv \min_{A_{L_0}, B_{L_0}} \cJ_0 + \Tr(\Lambda_0^{\top}\Lambda_0\Theta). \label{eqn:regularizer2} 
\end{align}
The Hessian matrix of the quadratic program (\textbf{P3}) is 
\begin{align*}
		H_{(\textbf{P3})} = \begin{bmatrix}
		\bar{x}_0\bar{x}_0^{\top} + \Theta & (\bar{x}_0\bar{x}_0^{\top}+ \Theta)C^{\top} \\ C(\bar{x}_0\bar{x}_0^{\top}+ \Theta) & C(\bar{x}_0\bar{x}_0^{\top}+ \Theta+X_0)C^{\top} + V
	\end{bmatrix}.
\end{align*}
Injecting $\theta_0$ into $x_0$ leads to a different ``regularization term'' as can be seen in \eqref{eqn:regularizer2}. By Lemma \ref{lemma:pdness}, it is immediate to check that (\textbf{P3}) is strictly convex, and hence has a unique minimum, and is also smooth. Subsequently, we know that the one-step filter $(A-L_0^*C, L_0^*)$ attains the minimum of $\cJ_0$; it also achieves a zero value for the regularization term $\Tr(\Lambda_0^{\top}\Lambda_0\Theta)$ in \eqref{eqn:regularizer2}, which is non-negative for any $\Theta >0$. Hence, $(A-L_0^*C, L_0^*)$ corresponds to the unique minimum of $(\textbf{P3})$. This proves that injecting $\theta_0$ into both $x_0$ and $\hat{x}_0$ results in the strictly convex quadratic problem $(\textbf{P3})$ whose unique minimum is exactly the one-step filter $(A-L_0^*C, L_0^*)$.

We proceed with the induction step. In the $\tau$-th iteration of the RHPG algorithm, for an arbitrary $\tau > 0$, assume that we have computed time-varying filters $\{A_{L_0}^*, \cdots, A_{L_{\tau-1}}^*, B_{L_{0}}^*, \cdots, B_{L_{\tau-1}}^*\}$. The quadratic program for the $\tau$-th iteration is constructed such that $\theta_0$ is injected into both $\hat{x}_{\tau}$ and $x_{\tau}$, but not $\hat{x}_0$ and $x_0$. We define
\begin{align*}
	\cJ_{\tau}:=\EE_{x_0, w_t, v_t} \Big\{ \sum_{t=0}^{\tau+1}(x_t-\hat{x}_t)^{\top}(x_t-\hat{x}_t)\Big\}.
\end{align*}
An associated quadratic program (\textbf{P}$_{\tau}$) can be represented as 
\begin{align*}
	(\textbf{P}_{\tau}) &\min_{A_{L_\tau}, B_{L_{\tau}}} \EE_{x_0, w_t, v_t, \theta_0} \Big\{ \sum_{t=0}^{\tau+1}(x_t-\hat{x}_t)^{\top}(x_t-\hat{x}_t)\Big\}, \\
	\text{s.t.} \quad &\hat{x}_{t+1} = A_{L_t}^*\hat{x}_t + B_{L_t}^*y_t, \quad \forall t \in \{0, \cdots, \tau-2\}, \quad \hat{x}_0 = \bar{x}_0\\
	& x_{t+1} = Ax_t + w_t, \quad \forall t \in \{0, \cdots, \tau-2\}, \ x_0 \sim \cN(\bar{x}_0, X_0)\\
	 &\hat{x}_\tau =  A_{L_{\tau-1}}^*\hat{x}_{\tau-1} + B_{L_{\tau-1}}^*y_{\tau-1} + \theta_0, \quad x_{\tau} = Ax_{\tau-1} + w_{\tau-1} + \theta_0.
\end{align*}
where the expectation in $\cJ_{\tau}$ is taken over $x_{\tau}$, $\hat{x}_{\tau}$, $w_{\tau}$, and $v_{\tau}$.
The Hessian matrix of the quadratic program (\textbf{P}$_\tau$) is 
\begin{align*}
		H_{(\textbf{P}_{\tau})} &= \EE_{x_{\tau}, \hat{x}_{\tau}, v_{\tau}}\begin{bmatrix}
		\hat{x}_{\tau}\hat{x}_{\tau}^{\top} & \hat{x}_{\tau}y_{\tau}^{\top} \\ y_{\tau}\hat{x}_{\tau}^{\top} & y_{\tau}y_{\tau}^{\top}
	\end{bmatrix} = \begin{bmatrix}
		\mu_{\hat{x}_{\tau}}\mu_{\hat{x}_{\tau}}^{\top} +\Theta & (\mu_{\hat{x}_{\tau}}\mu_{x_{\tau}}^{\top} +\Theta)C^{\top} \vspace{0.2em}\\ C(\mu_{x_{\tau}}\mu_{\hat{x}_{\tau}}^{\top} +\Theta) & C(\mu_{x_{\tau}}\mu_{x_{\tau}}^{\top} +\Theta)C^{\top} + V
	\end{bmatrix},
\end{align*}
where we have used $\mu_{\hat{x}_{\tau}}$ and $\mu_{x_\tau}$ to denote $\EE[\hat{x}_\tau]$ and $\EE[x_\tau]$, respectively. By Lemma \ref{lemma:pdness}, the Hessian matrix $H_{(\textbf{P}_{\tau})}$ is pd for any $\tau >0$, which proves that $(\textbf{P}_{\tau})$ is strictly convex and smooth. Furthermore, we can apply a derivation similar to \eqref{eqn:regularizer2} to show that $(\textbf{P}_{\tau})$ is equivalent to
\begin{align*}
	\min_{A_{L_\tau}, B_{L_\tau}} \cJ_{\tau} + \Tr(\Lambda_\tau^{\top}\Lambda_\tau\Theta), \quad  \Lambda_{\tau} := A-B_{L_\tau}C-A_{L_\tau}.
\end{align*}
The unique minimum to $(\textbf{P}_{\tau})$ is, therefore, $(A-L^*_{\tau}C, L^*_{\tau})$, which is the time-varying filter at $t=\tau$. By induction, we have proved Theorem \ref{theorem:quadratic}.

\end{document}